\documentclass[12pt]{article}

\usepackage{amssymb}    % use all AMS fonts
\usepackage{amsthm}     % AMS style theorem and proof environments

\newtheorem{theorem}{Theorem}[section]
\newtheorem{lemma}[theorem]{Lemma}
\newtheorem{proposition}[theorem]{Proposition}
\newtheorem{corollary}[theorem]{Corollary}

\newenvironment{prf}[1]{\trivlist
\item[\hskip
\labelsep{\it #1.\hspace*{.3em}}]}{%~\hspace{\fill}~$\square$%
\endtrivlist}

\newtheorem{predefinition}[theorem]{Definition}
\newenvironment{definition}{\begin{predefinition}\rm}{\end{predefinition}}
\newtheorem{preremark}[theorem]{Remark}
\newenvironment{remark}{\begin{preremark}\rm}{\end{preremark}}
\newtheorem{prenotation}[theorem]{Notation}

\newtheorem{preexample}[theorem]{Example}
\newenvironment{example}{\begin{preexample}\rm}{\end{preexample}}
\newtheorem{preclaim}[theorem]{Claim}

\newtheorem{prequestion}[theorem]{Question}

 \makeatletter
\def\emppsubsection{\@startsection{subsection}{2}{\z@}{-3.25ex plus -1ex minus -.2ex}{-1em}{\bf}}

\makeatother

\newcommand \C {{\cal C}}
\newcommand \X {{\cal X}}
\newcommand \Y {{\cal Y}}
\newcommand \B {{\cal B}}
\newcommand \CO {{\cal O}}

\newcommand \CM {{\cal M}}

\newcommand \ZZ {{\mathbb Z}}
\newcommand \NN {{\mathbb N}}

\newcommand  \FF {{\mathbb F}}
\newcommand  \RR {{\mathbb R}}

\newcommand  \QQ {{\mathbb Q}}

\newcommand \GG {{\mathbb G}}

\newcommand \Aut {\mathop{\rm Aut}}

\newcommand \Gal {\mathop{\rm Gal}}

\newcommand \gc {\mathop{\rm gcd}}
\newcommand \Hom {\mathop{\rm Hom}}
\newcommand \Ind {\mathop{\rm Ind}}

\newcommand \Spec {\mathop{\rm Spec}}

\newcommand \val {\mathop{\rm val}}

\title{Equiramified deformations of covers in positive characteristic}
\author{Rachel Pries: Colorado State University
\footnote{The author was partially supported by NSF grant DMS-04-00461.}}
\date{}

\begin{document}
\maketitle

\begin{abstract}
\noindent
Suppose $\phi$ is a wildly ramified cover of germs of curves 
defined over an algebraically closed field of characteristic $p$.
We study unobstructed deformations of $\phi$ in equal characteristic, which are equiramified 
in that the branch locus is constant and the ramification filtration is fixed.  
We show that the moduli space $\CM_\phi$ parametrizing equiramified deformations of $\phi$
is a subscheme of an explicitly constructed scheme.
This allows us to give an explicit upper and lower bound for the Krull dimension $d_\phi$ of $\CM_\phi$.
These bounds depend only on the ramification filtration of $\phi$.
When $\phi$ is an abelian $p$-group cover, 
we use class field theory to show that the upper bound for $d_\phi$ is realized.

\end{abstract}

\section{Introduction}

There are many open questions about Galois covers of curves in characteristic 
$p$ whose characteristic 0 analogues are well-understood.  
For example, consider a $G$-Galois cover $\varphi:\Y \to \X$ of Riemann surfaces. 
Deformations of $\varphi$ are parametrized by a moduli space
whose dimension, $3g_{\X}-3+|\B| +{\rm dim}\Aut(\X-\B)$, is determined by the genus 
$g_{\X}$ of $\X$ and the size $|\B|$ of the branch locus of $\varphi$.   
By the Riemann-Hurwitz formula, the genus of $\Y$ is determined by $|G|$, $g_{\X}$, $|\B|$ and the 
orders of the inertia groups.  

These statements are no longer true when $\varphi:\Y \to \X$ is a wildly ramified cover of curves.  
In characteristic $p$, the number of unramified covers of a fixed 
affine curve with a fixed Galois group is typically infinite.  
Not only can these covers be deformed without varying $\X$ or the branch
locus $\B$ of $\varphi$, but they can often be distinguished from each other by 
studying finer ramification invariants such as the conductor.
The genus of $\Y$ now depends on these finer ramification invariants. 

Let $k$ be an algebraically closed field of characteristic $p > 0$.
Suppose $\X$ is a smooth projective $k$-curve and $\B$ is a finite set of points of $\X$.
Suppose $G$ is a finite quotient of $\pi_1(\X-\B)$. 
(When $|\B|$ is nonempty, these groups have been classified by 
Raynaud \cite{Ra:ab} and Harbater \cite{Ha:ab} in their proof of Abhyankar's Conjecture).
Suppose $\varphi:\Y \to \X$ is a $G$-Galois cover of curves branched only at $\B$.

An important problem is to understand the deformation theory of $\varphi$.
It is well-known by the theory of formal patching \cite{Ha:long} 
that deformations of $\varphi$ satisfy a local/global property.
This theory allows one to simplify the question of deformations of $\varphi$
to the analogous question of deformations of an $I$-Galois cover $\phi:Y \to X$ of germs of curves.
Here $X$ is the germ of $\X$ at a branch point $b$, 
$Y$ is the germ of $\Y$ at a ramification point $\eta\in \varphi^{-1}(b)$
and $I$ is the inertia group of $\varphi$ at $\eta$. 
(Likewise, deformations of $\Y$ respecting the associated group action $G \to \Aut(\Y)$ 
satisfy a local/global property, \cite{BM}.)

\paragraph{Results.}

In this paper, we consider equal characteristic deformations of an $I$-Galois cover $\phi$ of germs of curves.
We study a functor $F_\phi$ parametrizing unobstructed deformations of $\phi$ which are equiramified 
in that the branch locus and ramification filtration do not change.  
The main result, Theorem \ref{Tsmoothdef2}, states that there is a moduli space $\CM_\phi$ 
representing the functor $F_\phi$ in a certain category.
Furthermore, $\CM_\phi$ is a subscheme of a direct product of schemes, 
each of which is a moduli space in its own right 
and which can be explicitly described in terms of the ramification filtration of $\phi$.   
This allows us to give an explicit upper and lower bound for the Krull dimension $d_\phi$ of $\CM_\phi$;
these bounds depend only on the ramification filtration of $\phi$. 

When $I$ is an abelian $p$-group, we strengthen this result in Corollary \ref{Tcyclic} using class field theory.  
In this case, $\CM_\phi$ is a direct product of copies of $\hat{\GG}_a$ modulo an action of $\FF_p^*$.
As a result, when $I$ is an abelian $p$-group, we give an exact formula for the Krull dimension $d_\phi$ of $\CM_\phi$ 
in terms of the breaks in the filtration of higher ramification groups of $\phi$.
Here the upper bound for $d_\phi$ is realized.
 
Corollary \ref{Creducible} shows that the upper bound for $d_\phi$ is also realized in the case that 
$I$ is a semi-direct product of the form $(\ZZ/p)^e \rtimes \ZZ/m$.
This type of inertia group occurs when the projective curve $\Y$ is ordinary.
When $\Y$ is ordinary, there is also a major restriction on the ramification filtration of $\phi$.
In this way, Corollary \ref{Creducible} is a generalization of 
part of the results of Cornelissen and Kato \cite{CK} who study deformations 
of the type arising from ordinary projective curves.

The upper bound for $d_\phi$ is not always realized. 
We give an example of this in Section \ref{Squaternion} when $p=2$ and when $I$ is the quaternion group of order 8. 

\paragraph{Outline.}

Section \ref{Sbackground} contains background information on the ramification 
of a wildly ramified cover $\phi$ of germs of curves.
We explain the connection between some deformations of $\phi$ and 
some $A \rtimes \mu_m$-covers where $A$ is an elementary abelian $p$-group.
These $A \rtimes \mu_m$-covers play an important role in later sections.
We give a detailed description of their equations in Section \ref{Seq}.

In Section \ref{Smoduli}, we construct a moduli space which parametrizes 
these $A \rtimes \mu_m$-covers in a certain category, Theorem \ref{Tmoduli}.  
It is a direct product of copies of $\hat{\GG}_a$ modulo an action by $\FF_{|A|}^*$.  
Its Krull dimension can be explicity computed in terms of the ramification filtration of the covers.
As an application, in Section \ref{CK}, we generalize the dimension count of \cite[Thm.\ 5.1(a)]{CK}. 

Sections \ref{Slocaldef} and \ref{Scategory} contain the precise definition of the deformation functor $F_\phi$.
In Section \ref{Ssubdef}, we use a Galois action from \cite{Pr:genus}
to study equiramified deformations of $\phi$ which have a constant quotient deformation.
This Galois action acts non-trivially only on an $A$-Galois subcover of $\phi$.
It causes the set of wildly ramified $I$-Galois covers which dominate a fixed quotient  
to form a principal homogeneous space 
under the action of a group of $A \rtimes \mu_m$-Galois covers.
The main results here are Theorem \ref{Tfunctor} and Corollary \ref{Cisofunc} which simplify the problem of 
equiramified deformations of $\phi$.  

In Section \ref{Stower}, we prove Theorem \ref{Tsmoothdef2} which states that there is a moduli space 
representing $F_\phi$ and gives an upper and lower bound for its Krull dimension $d_\phi$. 
The method is to study equiramified deformations of a tower of covers using induction.
The results from Section \ref{S2} allow one to reduce the question of deformations of $\phi$
to the question of deformations of each of the steps in the tower.  
However, it is difficult to verify whether the deformation is equiramified in a step-to-step manner.
We rephrase this issue in terms of an equiramified embedding problem.

The equiramified embedding problem has a solution in the case of abelian $p$-groups.
This yields the exact formula for $d_\phi$ in Corollary \ref{Tcyclic}.
We show in Section \ref{Squaternion} that the equiramified embedding problem does not always have a solution 
and thus the upper bound for $d_\phi$ is not always realized.

For a future application, we explain in Section \ref{BM} how the formula for $d_\phi$ can be used to study the 
Krull dimension of $\CM_g[G] \otimes \FF_p$.  This can be seen for $G=\ZZ/p$ in \cite[Section 5.2]{BM}.

I would like to thank D. Harbater and M. Raynaud for their invaluable feedback, 
along with J. Achter, I. Bouw, G. Cornelissen, A. Tamagawa and 
the participants of the conferences in Banff and in Leiden
for their helpful suggestions.  

\section{Galois covers of germs of curves} \label{Sbackground}

\subsection{Structure of the inertia group} \label{S1.1}

This section contains background material on the inertia group and ramification filtration
of a wildly ramified Galois cover of germs of curves.

Let $k$ be an algebraically closed field of characteristic $p >0$.
Consider an irreducible $k$-scheme $\Omega$.
Let $U=\underline{\Spec}(\CO_\Omega[[u]])$ 
and let $\xi$ be the closed point of $U$ defined by the equation $u=0$.

\paragraph{Inertia group.}

Suppose ${\phi}:{Y} \to U$ is a Galois cover of normal irreducible germs 
of $\Omega$-curves which is wildly ramified at the closed point $\eta ={\phi}^{-1}(\xi) \in {Y}$.
By \cite[Lemma 2.1.4]{Pr:deg}, after an \'etale pullback of $\Omega$, 
the decomposition group and inertia group over the generic point of $\eta$ are the same
and so the Galois group of ${\phi}$ is the same as its inertia group.  

Recall that the inertia group $I$ at the generic point of $\eta$ is of the form $P \rtimes_\iota \ZZ/m$ 
where $|P|=p^e$ for some $e > 0$ and $p \nmid m$, \cite[IV, Cor.\ 4]{Se:lf}. 
Here $\iota$ denotes the automorphism of $P$ which determines the 
conjugation action of $\ZZ/m$ on $P$.

\paragraph{Ramification filtrations.}

Associated to the cover ${\phi}$, there are two filtrations of $I$, namely 
the filtration of higher ramification groups $I_{\tilde{c}}$ in the lower numbering
and the filtration of higher ramification groups $I^c$ in the upper numbering. 
Let $\pi$ be a uniformizer of $\CO_{Y}$ at the generic point of $\eta$. 
If $\tilde{c} \in \NN$, then $I_{\tilde{c}}$ is the normal subgroup of all $g \in I$ 
such that $g$ acts trivially on $\CO_{Y}/\pi^{\tilde{c}+1}$.
Equivalently, $I_{\tilde{c}}=\{g \in I| \val(g(\pi)-\pi) \geq \tilde{c}+1\}$, \cite[IV, Lemma 1]{Se:lf}. 
If $\tilde{c} \in \RR^+$ and $c'=\lfloor \tilde{c} \rfloor$, then $I_{\tilde{c}}=I_{c'}$.   
Recall by Herbrand's formula \cite[IV, Section 3]{Se:lf}, 
that the filtration $I^c$ in the upper numbering is given by
$I^c=I_{\tilde{c}}$ where $\tilde{c}=\Psi(c)$ and $\Psi(c)=\int_0^c(I^0:I^t)dt$.  
Equivalently, $c=\int_0^{\tilde{c}}dt/(I_0:I_t)$.

We say that $j \in \NN^+$ is a {\it lower jump} of ${\phi}$ at $\eta$ if $I_{j} \not = I_{j+1}$.
In this case, $I_j/I_{j+1} \simeq (\ZZ/p)^{\ell_j}$ for some $\ell_j \in \NN^+$, \cite[IV, Cor.\ 3]{Se:lf}.
Then $\ell_j$ is the {\it multiplicity} of $j$.
A rational number $c$ is an {\it upper jump}
of ${\phi}$ at $\eta$ if $\Psi(c)=j$ for some lower jump $j$.   
Let $j_1, \ldots, j_e$ (resp.\ $\sigma_1, \ldots, \sigma_e$) be
the set of lower (resp.\ upper) jumps of ${\phi}$ at $\eta$
written in increasing order with multiplicity.  
These are the positive breaks in the filtration of ramification groups in the lower 
(resp.\ upper) numbering.
By \cite[IV, Prop.\ 11]{Se:lf}, $p \nmid j_i$ for any lower jump $j_i$.
Herbrand's formula implies that $j_i-j_{i-1}=(\sigma_i-\sigma_{i-1})|I|/|I_{j_i}|$. 
%In particular, $\sigma_i|I|/|I^{\sigma_i}| \in \NN$ and $\sigma_1=j_1/m$.

The number $\sigma=\sigma_e$ is the {\it conductor} of 
${\phi}$ at $\eta$; $\sigma$ is the largest $c \in \QQ$ such that 
inertia group $I^c$ is non-trivial in the filtration of higher ramification groups in the upper numbering.  
(This indexing is slightly different
than in \cite{Se:lf}, where the ideal $(x^{\sigma+1})$ is the {\it conductor} of the 
extension of complete discrete valuation rings for a given uniformizer $x$ at the branch point.) 
 
The ramification filtration in the upper numbering, the upper jumps, and the conductor
are preserved under quotients, \cite[IV, Prop.\ 14]{Se:lf}.

\paragraph{An initial step in the filtration.}

We will study the deformations of ${\phi}$ by first studying the ones
which fix its $I/A$-Galois quotient for a suitable choice of $A \subset I$.

\begin{lemma} \label{LhypA}  Suppose ${\phi}:{Y} \to U$ is an $I$-Galois cover 
with inertia group $I = P \rtimes_\iota \ZZ/m$ and conductor $\sigma$ as above.
Then there exists $A \subset I^\sigma$ satisfying the following conditions: 
$A$ is central in $P$; $A$ is normal in $I$; 
$A$ is a nontrivial elementary abelian $p$-group; and $A$ is irreducible under the action of $\ZZ/m$.
\end{lemma}

\begin{proof}
In \cite[Lemma 2]{Pr:genus}, we show this follows immediately from \cite[IV]{Se:lf}.
\end{proof}

We fix $A \subset P$ satisfying the conditions of Lemma \ref{LhypA}.
Let $a$ be the positive integer such that $A \simeq (\ZZ/p)^a$ and let $q=p^a$.
Let $A \rtimes_\iota \ZZ/m$ be the semi-direct product determined by the 
restriction of the conjugation action of $\ZZ/m$ on $P$.
Let $\overline{I}=I/A$ and $\overline{P}=P/A$.

\subsection{Structure of the Galois cover}

\paragraph{Factoring the Galois cover.}

Consider an $I$-Galois cover ${\phi}:Y \to U$ of germs of $\Omega$-curves with conductor $\sigma$
and a subgroup $A$ satisfying the conditions of Lemma \ref{LhypA}.
This situation yields a factorization of ${\phi}$ which we denote 
$$Y \stackrel{\phi^A}{\to} \overline{Y} \stackrel{\overline{\phi}}{\to} X 
\stackrel{\kappa}{\to} U.$$
Here ${\phi^A}$ is $A$-Galois and $\overline{\phi}$ is $\overline{P}$-Galois. 
Also $\kappa \circ \overline{\phi}$ is $\overline{I}$-Galois.    
Let $\phi^P = \overline{\phi} \circ \phi^A$ denote the $P$-Galois subcover $Y \to X$ of ${\phi}$. 

The cover $\kappa:X \to U$ is a $\ZZ/m$-Galois cover.
By Kummer theory, $X \simeq \underline{\Spec}(\CO_\Omega[[x]])$ for some $x$ such that $x^m=u$.
The generator of $\ZZ/m$ takes $x \mapsto \zeta(x)$ for some $\zeta \in k$ with order $m$.
We choose a compatible system of roots of unity of $k$ so that $\zeta=\zeta_m$.
In this way, $\kappa$ yields an isomorphism $\ZZ/m \simeq \mu_m$. 
Let $\xi'=\kappa^{-1}(\xi) \in X$ be the ramification point of $\kappa$ and the branch point of $\phi^P$. 

Let $K$ (resp.\ $\overline{K}$) be the function field of $Y$ (resp.\ $\overline{Y}$).  
By \cite[Prop.\ 1.1]{GS}, there exists $y \in K$ and $r_\phi \in \overline{K}$ so that
the equation for $\phi^A$ is $y^q-y=r_\phi$.

\paragraph{Explanation of a Galois action on covers.} \label{S1.2}

The main point of Section \ref{Ssubdef} is that deformations of an $I$-Galois cover with constant $I/A$-Galois quotient 
are related to $A$-Galois covers of a certain type.  This relationship relies heavily on
a Galois action on wildly ramified germs of curves from \cite[Section 2]{Pr:genus}.
Before making this relationship more precise, we give a brief explanation of the Galois action.

Naively speaking, one can try to modify $\phi$ by changing $\phi^A$ while keeping $\overline{\phi}$ constant.  
This can be done by modifying the equation for $\phi^A$ to an equation of the form $y^q-y=r_\phi + r_\alpha$
for some $r_\alpha \in \overline{K}$.
By \cite[Lemma 3]{Pr:genus}, the composition $\overline{\phi} \circ \phi^A$ will still be $P$-Galois if and only if 
$r_\alpha$ is in the function field $\CO_\Omega[[x]][x^{-1}]$ of $X$;
in addition, every $P$-Galois cover dominating $\overline{\phi}$ occurs for some choice of $r_\alpha$ in the function field of $X$.
In this way, any modification of $\phi$ with fixed $I/A$-quotient corresponds to a 
(possibly disconnected) $A$-Galois cover $\psi:V \to X$ with equation $v^q-v=r_\alpha$.

Furthermore, by \cite[Lemma 5]{Pr:genus}, the composition $\kappa \circ \overline{\phi} \circ \phi^A$
will still be $I$-Galois if and only if $\kappa \circ \psi$ is an $A \rtimes_\iota \mu_m$-Galois cover.
Viewing $r_\alpha$ as a Laurent series in $x$, this places constraints on which of its coefficients can be non-zero.
The ramification filtration of $\phi$ will not change under this modification
exactly when the conductor of $\psi$ is bounded above by $m\sigma$, by \cite[Proposition 7]{Pr:genus}.

To summarize, this explains the basic idea of Section \ref{Ssubdef}: deformations of an $I$-Galois cover of curves
with constant $I/A$-Galois quotient correspond to $A$-Galois covers $\psi:V \to X$ of germs of curves 
having conductor at most $m\sigma$ so that $\kappa \circ \psi$ is an $A \rtimes_\iota \mu_m$-Galois cover. 
Luckily, it is easy to describe such $A \rtimes_\iota \mu_m$-Galois covers.

\subsection{Equations for $A \rtimes_\iota \mu_m$-Galois covers} \label{Seq}

Suppose $\psi:V \to X$ is an $A$-Galois cover of germs of $\Omega$-curves 
so that $\kappa \circ \psi$ is an $A \rtimes_\iota \mu_m$-Galois cover.
By \cite[Prop.\ 1.1]{GS}, there exists $v$ in the function field of $V$ so that the equation for $\psi$ is $v^q-v=r(x)$.
Here $r(x) \in \CO_\Omega[[x]][x^{-1}]$.

The conductor of $\psi$ can be determined explicitly in terms of $r(x)$.
The Laurent series $r(x)$ can be written as a sum of terms $(r_t)^{p^t}$ for $t \in \NN$ 
where $r_t \in \CO_\Omega[[x]][x^{-1}]$ and where $p$ divides no exponent of $x$ in $r_t$. 
Let $s_t$ be the degree of $r_t$ in $x^{-1}$.
We define the {\it prime-to-$p$ degree} of $r(x)$ to be $-{\rm min}\{s_t\ | \ t \in \NN\}$.
For example, 3 is the prime-to-$p$ degree of $x^{-2} + x^{-3p^2}$.

\begin{lemma} \label{Lppdegree}
There is exactly one break in the filtration of ramification groups of $\psi$.
The conductor $s$ of $\psi$ equals the prime-to-$p$ degree of $r(x)$.
\end{lemma}

Recall that $s \in \NN^+$ by the Hasse-Arf Theorem \cite[V, Thm.\ 1]{Se:lf} and that $p \nmid s$.

\begin{proof}
The first claim follows from the fact that $\mu_m$ acts irreducibly on $A$. 
For the second claim, by \cite[Prop.\ 1.2]{GS}, 
there is a $\ZZ/p$-Galois cover $\psi':V' \to X$ which is a quotient of $\psi$ and has equation $v_1^p-v_1=r(x)$. 
The conductor of $\psi$ equals the conductor of $\psi'$. 
The conductor and lower jump of $\psi'$ are the same by Herbrand's formula. 
There is an isomorphism (over the generic geometric fibre of $\Omega$) between $\psi'$ and the $\ZZ/p$-Galois cover 
$v_2^p-v_2=r(x)'$ where $r(x)'= \sum_{t \in \NN} r_{t}$ by Artin-Schreier theory.
By definition, this degree is prime-to-$p$ and equals the prime-to-$p$ degree of $r(x)$. 
By \cite[Prop.\ VI.4.1]{St}, the lower jump of $\psi'$ is the degree of $r(x)'$ in $x^{-1}$.
\end{proof}

The congruence value of the conductor modulo $m$ can be determined from the 
semi-direct product $A \rtimes_\iota \mu_m$. 

\begin{lemma} \label{Lcong}
Associated to the $\mu_m$-Galois cover $\kappa: X \to U$ and the group $A \rtimes_\iota \mu_m$,
there is a unique integer $s_\iota$ (such that $1 \leq s_\iota \leq m$) with the following property:
for any $A$-Galois cover $\psi:V \to X$ with conductor $s$ so that 
$\kappa \circ \psi$ is an $A \rtimes_\iota \mu_m$-Galois cover then $s \equiv s_\iota \bmod m$. 
Furthermore: 
\begin{description}
\item{(i)} $\gc(m,s_\iota)=|{\rm Ker}(\iota)|$;

\item{(ii)} $[\FF_p(\zeta_m^{s_\iota}):\FF_p]=a$;

\item{(iii)} If $cx^i$ is a term of $r(x) \in \CO_\Omega[[x]][x^{-1}]]$ for which $c \not = 0$,
then $i \equiv s_\iota \bmod m$.
\end{description}
\end{lemma}

\begin{proof}
Suppose $\psi:V \to X$ is an $A$-Galois cover so that 
$\kappa \circ \psi$ is an $A \rtimes_\iota \mu_m$-Galois cover.
By the definition of $A$-Galois cover, there is an isomorphism between $A$ and $\Gal(V/X)$. 
There is a natural isomorphism between $\Gal(V/X)$ and $\FF_q$
where $\gamma(v)=v+\gamma$ for $\gamma \in \FF_q$.   
By the definition of $\kappa$, there is an isomorphism between $\mu_m$ and $\Gal(X/U)$
where the generator takes $x \mapsto \zeta_m x$. 

The semi-direct product $A \rtimes_\iota \mu_m$ determines an automorphism $\iota(\zeta_m) \in \Aut(A)$ so that
$\iota(\zeta_m)(\tau)= \zeta_m \tau \zeta_m^{-1}$ for $\tau \in A$. 
By \cite[IV, Prop.\ 6.a]{Se:lf}, $\zeta_m \tau \zeta_m^{-1}=z \tau$ for some $z \in \mu_m$.
Let $s_\iota$ be the integer such that $\zeta_m^{s_\iota}=z$ and $1 \leq s_\iota \leq m$. 

Suppose $s$ is the conductor of $\psi$. 
For any $\tau \in A$, the lower jump of $\psi$ for $\tau$ is equal to $s$ by Lemma \ref{Lppdegree} and Herbrand's Theorem. 
Furthermore, $\zeta_m^s=z$ by \cite[IV, Prop.\ 9]{Se:lf}.
Thus $s \equiv s_\iota$ modulo $m$. 
\begin{description}
\item{(i)} For $\zeta \in \mu_m$, $\iota(\zeta)(\tau)=\zeta^{s_\iota} \tau$. 
Thus ${\rm Ker}(\iota)$ is the subgroup of $\mu_m$ consisting of elements whose order divides $s_\iota$.
This subgroup has size $\gc(m, s_\iota)$.

\item{(ii)} First, $\zeta_m^{s_\iota} \in \FF_q$ since $\iota$ stabilizes $A \simeq \FF_q$.
Second, the fact that $A$ is irreducible under $\iota$ implies that 
$\zeta_m^{s_\iota} \not \in \FF_{p^{a'}}$ for any $a' <a$.  
Thus $\FF_p(\zeta_m^{s_\iota})=\FF_q$.

\item{(iii)} The equation $v^q-v=r(x)$ is invariant under the action $\zeta_m: x \mapsto \zeta_m x$.
This implies that $\zeta_m:v \mapsto \zeta_m^{-s_\iota} v$ and that the exponents of $r(x)$ are all 
congruent modulo $m$. 
\end{description}
\end{proof}

\begin{lemma} \label{Lnoniso} 
For $i=1,2$, suppose $\psi_i:V_i \to X$ is an $A$-Galois cover with equation $v_i^q-v_i=r_i(x)$
where $r_i(x) \in \CO_\Omega[[x]][x^{-1}]$. 
There is an isomorphism between the $A$-Galois covers $\psi_1$ and $\psi_2$ if and only if 
$r_2(x)=\zeta r_1(x) +d^q-d$ for some $\zeta \in \FF_q^*$ and some $d \in \CO_\Omega[[x]][x^{-1}]$.
\end{lemma}

One can also show that $\kappa \circ \psi_1:V_1 \to U$ and $\kappa \circ \psi_2:V_2 \to U$ 
are isomorphic as $A \rtimes_\iota \mu_m$-Galois covers if and only if $\psi_1$ and $\psi_2$ are isomorphic as $A$-Galois covers.

\begin{proof}
Suppose $r_1(x)$ and $r_2(x)$ are such that there is an isomorphism $\beta$ between the $A$-Galois covers $\psi_1$ and $\psi_2$.
Since $\beta$ is invertible, $\beta(v_2)=\zeta v_1 +d$ for some unit $\zeta$ 
and some $d$ in the function field $\CO_\Omega[[x]][x^{-1}]$ of $X$. 
Thus $\zeta^q v_1^q -\zeta v_1 +d^q- d= r_2(x)$.
Thus $(\zeta^q-\zeta)v_1= - \zeta^qr_1(x) -(d^q-d)+ r_2(x)$.
This implies that $(\zeta^q-\zeta)v_1$ is contained in the function field of $X$ 
which is only possible if $\zeta \in \FF_q^*$.  
Then $r_2(x)= \zeta r_1(x) + d^q-d$.

For the converse, if $r_2(x)=\zeta r_1(x) +d^q-d$, 
then there is an isomorphism $\beta$ between $\psi_1$ and $\psi_2$ as $A$-Galois covers where
$\beta(v_2)=\zeta v_1 +d$.
\end{proof}
 
An equation $v^q-v=r(x)$ for an $A$-Galois cover $\psi:V \to X$ 
is in {\it standard form} if $r(x) \in \CO_\Omega[x^{-1}]$ and if
$q$ does not divide the exponent of any term of $r(x)$.
Lemma \ref{Lcong} implies that the congruence value of the exponents of $r(x)$ modulo $m$
does not change in the process of changing the equation for $\psi$ into standard form.
To see this, suppose $cx^i$ is a term of $r(x)$ for which $i=qi_1$ and $c \not = 0$.
Then (iii) implies that $\gc(m,s_\iota)$ divides $i$ and thus divides $i_1$.
Also (ii) implies that $m/\gc(m,s_\iota)$ divides $q-1$.  Thus $i \equiv i_1 \bmod m$.

\section{Deformations via a Galois action} \label{S2}

We study the deformations of an $I$-Galois cover $\phi_k: Y_k \to U_k$ of germs of $k$-curves
whose $I/A$-Galois quotient is constant.  
The equations for these deformations can be related via a Galois action to 
the equations for $A$-Galois covers which are easy to understand.
In fact, we prove that the functor $F_{\phi, A}$ which parametrizes these deformations is 
isomorphic to another functor which parametrizes the $A$-Galois covers.
We give an explicit description of the moduli space that represents the latter functor 
in a certain category and compute its Krull dimension.

\subsection{Definition of deformation functor} \label{Slocaldef}

Suppose $\phi=\phi_k:Y_k \to U_k$ is an $I$-Galois cover of normal irreducible germs of
$k$-curves which is wildly ramified at the closed point $\eta_k = \phi^{-1}(\xi_k)$.
Here $U_k=\Spec(k[[u]])$ is the germ of a smooth $k$-curve at a closed point $\xi_k$.

Let $\C$ be the category of irreducible pointed $k$-schemes; we denote an object of $\C$ by 
$(\Omega, \omega)$ where $\Omega$ is an irreducible $k$-scheme and $\omega$ is the image of the 
chosen morphism $\Spec(k) \to \Omega$.
Let $U_\Omega=\underline{\Spec}(\CO_\Omega[[u]])$  
and let $\xi_\Omega$ be the closed point defined by the equation $u=0$.
We denote by $U_\omega$ the fibre of $U_\Omega$ over $\omega$.
For each object $(\Omega, \omega)$ of $\C$, there is a natural
morphism $\CO_\Omega[[u]] \to k[[u]]$ taking $u \mapsto u$.
This yields a natural choice of isomorphism $U_k \stackrel{\sim}{\to} U_\omega$.

A {\it deformation} of $\phi$ over $(\Omega, \omega)$ is an $I$-Galois cover
$\phi_{\Omega}: Y_{\Omega} \to  U_\Omega$ of normal irreducible germs of $\Omega$-curves
together with an isomorphism between $\phi$ and the fibre of $\phi_{\Omega}$ over $\omega$ 
as $I$-Galois covers of $U_\omega$.
Two deformations $\phi_{\Omega}$ and $\phi_{\Omega}'$ over $(\Omega, \omega)$ are 
isomorphic if there is an isomorphism between them
as $I$-Galois covers of the $\Omega$-curve $ U_\Omega$
which commutes with the isomorphisms between $\phi$ and their fibres over $\omega$.
A deformation is {\it constant} if it is isomorphic to the cover of germs of $\Omega$-curves with 
constant fibres, namely the restriction over $U_\Omega$
of $\phi \times_k \Omega:Y_k \times_k \Omega \to U_k \times_k \Omega$.

A deformation $\phi_{\Omega}$ is {\it equiramified} if
its branch locus consists of only the $\Omega$-point $\xi_{\Omega}$
and if it is totally ramified over the generic point of $\xi_{\Omega}$. 
For an equiramified deformation, a result of Kato \cite{Ka} implies
that the degree of the ramification divisor is the same for the fibre over $\omega$ 
and for the generic geometric fibre, \cite[Lemma 11]{Pr:genus}.

\begin{definition}
Let $F_{\phi}$ (resp.\ $F_{\phi, A}$) be the contravariant deformation functor from $\C$ to sets
which associates to $(\Omega,\omega)$ the set of isomorphism classes of equiramified deformations of
$\phi$ over $(\Omega,\omega)$ (resp.\ for which the $I/A$-Galois quotient deformation is constant).
\end{definition} 

The main result of Section \ref{Ssubdef} is that the functor $F_{\phi, A}$ 
is isomorphic to another functor $F_{A \rtimes_\iota \mu_m, \sigma}$ which we define now.

For an irreducible pointed $k$-scheme $(\Omega, \omega)$,
let $X_\Omega=\underline{\Spec}(\CO_\Omega[[x]])$.
Consider the $\mu_m$-Galois cover $\kappa:X_\Omega \to U_\Omega$ of germs of $\Omega$-curves with equation $x^m=u$  
and Galois action $x \mapsto \zeta_m x$ as in Section \ref{S1.2}.
Let $\xi'_\Omega$ be the closed point defined by the equation $x=0$ and let $X'_\Omega=X_\Omega-\{\xi'_\Omega\}$.

Consider the group $H_A(\Omega)=\Hom(\pi_1(X'_\Omega), A)$. 
We suppress the choice of basepoint of $X'_\Omega$ from the notation.
An element $\alpha \in H_A(\Omega)$ may be identified with the isomorphism class of 
an $A$-Galois cover $\psi_\alpha:V \to X_\Omega$ branched only over the closed point $\xi'_\Omega$.
The identity $\alpha_{0, \Omega}$ of $H_A(\Omega)$ corresponds to the totally disconnected $A$-Galois cover
$\psi_{0,\Omega}:\Ind_{0}^A X_\Omega \to X_\Omega$.

The automorphism $\iota$ of $A$ induces a natural automorphism $\overline{\iota}$
of $H_A(\Omega)$.
Let $H_A^\iota(\Omega)$ be the subgroup of $H_A(\Omega)$ fixed by $\overline{\iota}$.
In other words, $\alpha$ is an element of $H_A^\iota(\Omega)$ if and only if
the composition $\kappa \circ \psi_\alpha:V \to U_\Omega$ is an $(A \rtimes_\iota \mu_m)$-Galois cover.

The {\it conductor} of $\alpha \in H_A(\Omega)$ is the conductor 
of $\psi_\alpha$ over its generic geometric fibre.
Given $\sigma \in \QQ^+$, let $H_A^{\iota, \sigma}(\Omega)$ be the subset of $H_A^\iota(\Omega)$
consisting of elements $\alpha$ for which the conductor of $\psi_\alpha$ is at most $m\sigma$.
In fact, $H_A^{\iota, \sigma}(\Omega)$ is a subgroup of $H_A^\iota(\Omega)$.  

\begin{definition}
Let $F_{A \rtimes_\iota \mu_m, \sigma}$ be the contravariant $A$-Galois functor from $\C$ to sets
which associates to $(\Omega,\omega)$ the subset of elements $\alpha \in H_A^{\iota, \sigma}(\Omega)$
for which the fibre of $\psi_\alpha$ over $\omega$ is isomorphic to the totally disconnected cover $\psi_{0, \omega}$.
\end{definition}

\subsection{Simplification of the deformation functor} \label{Ssubdef}

We show that the functor $F_{\phi, A}$ parametrizing deformations of the $I$-Galois cover $\phi$
with constant $I/A$-Galois quotient is isomorphic to the functor $F_{A \rtimes_\iota \mu_m, \sigma}$ 
parametrizing $A$-Galois covers of a certain type.
The advantage of this is that we can explicitly describe a moduli space for 
the functor $F_{A \rtimes_\iota \mu_m, \sigma}$ in a certain category, Section \ref{Smoduli}.
   
\begin{theorem} \label{Tfunctor}
The functors $F_{\phi, A}$ and $F_{A \rtimes_\iota \mu_m, \sigma}$ are isomorphic.
\end{theorem}

Let $(\Omega, \omega)$ be a pointed $k$-scheme.  
For the proof, we show that there is a functorial bijection between the sets
$F_{\phi, A}(\Omega, \omega)$ and $F_{A \rtimes_\iota \mu_m, \sigma}(\Omega, \omega)$.
In fact, the proof shows that $F_{\phi, A}(\Omega, \omega)$ and $F_{A \rtimes_\iota \mu_m, \sigma}(\Omega, \omega)$ 
are isomorphic as groups, where the identity element of $F_{\phi, A}(\Omega, \omega)$ 
is the constant deformation of $\phi$.  

\begin{proof}
Let $\phi^P_{0,\Omega}$ (resp.\ $\overline{\phi}_{0,\Omega}$) denote the 
constant deformation of the $P$-Galois subcover $\phi^P$ 
(resp.\ of the $\overline{P}$-Galois subquotient $\overline{\phi}$) of $\phi$ over $\Omega$.

The reduction of a $P$-Galois cover modulo $A$ yields a $\overline{P}$-Galois cover.
Consider the resulting morphism 
${\rm red}_A: \Hom(\pi_1(X'_\Omega), P) \to \Hom(\pi_1(X'_\Omega), \overline{P})$.
Consider the fibre $H_{\overline{\phi}}(\Omega)$ of ${\rm red}_A$ over $\overline{\phi}_{0,\Omega}$.
Then $\phi^P_{0,\Omega}$ corresponds to an element of this fibre.

Let $H_{A}(\Omega)=\Hom(\pi_1(X'_\Omega), A)$.
By \cite[Lemma 3]{Pr:genus}, $H_A(\Omega)$ acts simply transitively on the fibre $H_{\overline{\phi}}(\Omega)$. 
Equivalently, the fibre $H_{\overline{\phi}}(\Omega)$ is a principal homogeneous space for $H_{A}(\Omega)$.
Given an element $\phi^P_\Omega$ of the fibre $H_{\overline{\phi}}(\Omega)$, 
it follows that there exists exactly one $\alpha \in H_A(\Omega)$ so that $\phi^P_\Omega$ is the image
of $\phi^P_{0,\Omega}$ under the action of $\alpha$.
Conversely, given any $\alpha \in H_A(\Omega)$, 
the image of $\phi^P_{0,\Omega}$ under the action of $\alpha$ 
yields a well-defined element $\phi^P_\Omega$ of the fibre $H_{\overline{\phi}}(\Omega)$.
It follows that there is a functorial bijection between $P$-Galois covers 
$\phi^P_\Omega$ of normal germs of curves which are \'etale over $X'_\Omega$ and which
dominate the $\overline{P}$-Galois cover $\overline{\phi}_{0,\Omega}$
and elements $\alpha \in H_A(\Omega)$.

Let $\phi_\Omega$ denote an element of $F_{\phi, A}(\Omega, \omega)$.
By definition, $\phi_\Omega$ is (the isomorphism class of) 
an equiramified deformation of $\phi$ with constant $I/A$-Galois quotient deformation.
In particular, the $P$-Galois subcover $\phi^P_\Omega$ of $\phi_\Omega$
is a $P$-Galois cover of normal germs of curves which is \'etale over $X'_\Omega$ and which
dominates the $\overline{P}$-Galois cover $\overline{\phi}_{0,\Omega}$.
As above, it yields an element $\alpha \in H_A(\Omega)$.
We now show that $\alpha \in F_{A \rtimes_\iota \mu_m, \sigma}(\Omega, \omega)$.

We first consider the invariance of these covers under the $\mu_m$-Galois action.  
The composition $\kappa \circ \phi^P_{0,\Omega}$ is an $I$-Galois cover.
By \cite[Lemma 5]{Pr:genus},
the cover $\kappa \circ \phi^P_{\Omega}$ is an $I$-Galois cover if and only if $\alpha \in H_A^\iota(\Omega)$. 
Thus $\alpha \in H_A^\iota(\Omega)$. 

Let $s$ be the conductor of $\alpha$.
The fact that $\phi_\Omega$ is an equiramified deformation of $\phi$ implies that 
(over the complete local ring of $\Omega$ at $\omega$)
it satisfies the hypotheses of \cite[Lemma 11]{Pr:genus}. 
As a result, the degree of the ramification divisor of $\phi_\Omega$ is the same
over the generic geometric point of $\Omega$ and over $\omega$. 
In particular, the conductor $\sigma_\Omega$ of $\phi_\Omega$
equals the conductor $\sigma$ of $\phi_\omega$.
On the other hand, since $\phi_\Omega$ is irreducible, its conductor $\sigma_\Omega$
is the maximum of the conductor of $\phi_{0,\Omega}$ (namely $\sigma$)
and the conductor of $\kappa \circ \alpha$ (namely $s/m$) by \cite[Proposition 7]{Pr:genus}.
Thus $s \leq m\sigma$ and $\alpha \in H_A^{\iota, \sigma}(\Omega)$. 

Since $\phi_\Omega$ is a deformation of $\phi$,
the fibre of $\phi_\Omega$ and of $\phi_{0,\Omega}$ over $\omega$ are isomorphic. 
This implies the fibre of $\psi_\alpha$ over $\omega$ is isomorphic to $\psi_{0,\omega}$
and so $\alpha \in F_{A \rtimes_\iota \mu_m, \sigma}(\Omega, \omega)$.

Conversely, an element $\alpha \in F_{A \rtimes_\iota \mu_m, \sigma}(\Omega, \omega)$
determines as above a $P$-Galois cover $\phi^P_\Omega$ of normal germs of curves 
dominating $\overline{\phi}_{0,\Omega}$.
Let $\phi_\Omega=\kappa \circ \phi^P_\Omega$. 
Then $\phi_\Omega$ is $I$-Galois by \cite[Lemma 5]{Pr:genus}.
The fact that the fibre of $\psi_\alpha$ over $\omega$ is isomorphic to $\psi_{0,\omega}$
implies that the fibre of $\phi_\Omega$ over $\omega$ is isomorphic to $\phi$.
Then $\phi_\Omega$ is a cover of irreducible germs of curves since $\phi$ has that property.
As a result, $\phi_\Omega$ is a deformation of $\phi$
with constant $I/A$-quotient deformation.
Also $\phi_\Omega$ is equiramified by \cite[Proposition 7]{Pr:genus}.
Thus the isomorphism class of $\phi_\Omega$ is an element of $F_{\phi, A}(\Omega, \omega)$.
\end{proof}

\subsection{Equivalence classes of deformations} \label{Scategory}

One can construct a moduli space of finite dimension representing the functor $F_{A \rtimes_\iota \mu_m, \sigma}$ 
(and thus also $F_{\phi,A}$)
but only at the expense of leaving the category of pointed $k$-schemes.
Here are two examples that illustrate how a moduli space 
for deformations of a wildly ramified cover depends critically on the choice of category.

\begin{example} Suppose $\phi$ is the $\ZZ/p$-Galois cover $y^p-y=u^{-j}$.  
Let $\Omega=\Spec(k[t])$ and let $\omega$ be the point $t=0$.  
The equation $y_t^p-y_t=u^{-j}+t$ gives a deformation $\phi_\Omega$ 
of $\phi$ over $(\Omega, \omega)$.
Let $\Omega' \to \Omega$ be the \'etale cover $z^p-z=t$.
The choice of a point $\omega'$ over $\omega$ corresponds to the choice of a root of $z^p-z$.
There is an isomorphism between the pullback of $\phi_\Omega$ to $(\Omega', \omega')$ 
and the cover with equation $y^p-y=u^{-j}$ (where this isomorphism identifies $y$ with $y_t-z$).
So this pullback is isomorphic to the constant deformation of $\phi$ over $(\Omega', \omega')$.   
\end{example}

This example illustrates how the presence of nontrivial \'etale covers in the category $\C$ 
of pointed $k$-schemes makes it difficult to construct a moduli space which is fine, not coarse.
For this reason, we switch to working with equal characteristic complete local rings.

Let $\hat{\C}$ be the category of spectra of equal characteristic complete local rings.
It is clear from Section \ref{Slocaldef} how to define the functors 
$F_\phi$, $F_{\phi, A}$, and $F_{A \rtimes_\iota \mu_m, \sigma}$ on $\hat{C}$.
In the category $\hat{\C}$ all deformations will automatically be unobstructed.  

\begin{example}
Suppose $\phi$ is the $\ZZ/p$-Galois cover $y^p-y=u^{-j}$.
In the category $\hat{C}$, 
one can construct a family of deformations of $\phi$ over $\Spec(k[[t]])$ with infinite dimension
using the equations $y^p-y=u^{-j}+ t u^{-p^i}$ for $i \in \NN$.
\end{example}

The problem is that these deformations are all related via purely inseparable extensions of $k[[t]]$.  
As in \cite[Theorem 2.2.10]{Pr:deg}, one can resolve this difficulty in the category $\hat{\C}$
by constructing a {\it configuration space} instead of a moduli space.
Here we instead choose to work with a category where finite purely inseparable morphisms are invertible. 
In this category, there is a moduli space of finite Krull dimension representing the functors $F_{\phi, A}$ and 
$F_{A \rtimes_\iota \mu_m, \sigma}$, Section \ref{Smoduli}.

\begin{definition} \label{Dcat}
Let $\C'$ be the category whose objects are the objects of the category $\hat{\C}$ 
and whose morphisms consist of all morphisms in $\hat{\C}$ along with formal inverses to 
finite purely inseparable morphisms between objects in $\hat{\C}$. 
\end{definition}

Note that any such finite purely inseparable morphism (sometimes called a radicial morphism) 
is a composition of Frobenius morphisms \cite[IV.2.5]{Har}. 
Thus the category $\C'$ can be obtained by localizing the category $\hat{\C}$ by 
the multiplicative system of morphisms
which are powers of Frobenius \cite[Prop.\ 3.1]{Har:RD}. 

Suppose $\Omega$ and $\Omega'$ are objects of $\hat{C}$.
Suppose $\phi_{\Omega}$ (resp.\  $\phi_{\Omega'}$) is a deformation of $\phi$ 
over $\Omega$ (resp.\ $\Omega'$).
The deformations $\phi_{\Omega}$ and $\phi_{\Omega'}$ are {\it equivalent} 
if there exists an object $\Omega''$ of $\hat{C}$ along with 
finite purely inseparable (possibly trivial) morphisms 
$\pi:\Omega'' \to \Omega$ and $\pi':\Omega'' \to \Omega'$
so that the pullbacks $\pi^* \phi_{\Omega}$ and $(\pi')^* \phi_{\Omega'}$ are isomorphic 
deformations of $\phi$ over $\Omega''$.  
In particular, if $i: \Omega' \to \Omega$ is a finite purely inseparable morphism 
then the pullback $i^* \phi_\Omega$ of $\phi_{\Omega}$ is a deformation of $\phi$
over $\Omega'$ which is equivalent to the deformation $\phi_\Omega$.

\begin{example}
Suppose $\phi$ is the $\ZZ/p$-Galois cover $y^p-y=u^{-j}$.
The deformations $y^p-y=u^{-j} + tu^{-p}$ over $\Spec(k[[t]])$ 
and $y^p-y=u^{-j}+su^{-1}$ over $\Spec(k[[s]])$ are equivalent.   
In this case, $\pi: \Spec(k[[s]]) \to \Spec(k[[t]])$ is given by $t \mapsto s^p$ and 
$\pi':  \Spec(k[[s]]) \to \Spec(k[[s]])$ is the identity.
\end{example}

Likewise, suppose $\alpha$ (resp.\ $\alpha'$) is an element of $F_{A \rtimes_\iota \mu_m, \sigma}(\Omega)$ 
(resp.\ $F_{A \rtimes_\iota \mu_m, \sigma}(\Omega')$) with associated cover $\psi_\alpha$ (resp.\  $\psi_\alpha'$).
We say that $\alpha$ and $\alpha'$ are {\it equivalent}  
if there exists an object $\Omega''$ of $\hat{C}$ along with 
finite purely inseparable (possibly trivial) morphisms 
$\pi:\Omega'' \to \Omega$ and $\pi':\Omega'' \to \Omega'$
so that the pullbacks $\pi^* \phi_{\Omega}$ and $(\pi')^* \phi_{\Omega'}$ are isomorphic 
covers of $X_{\Omega''}$.  

Consider the contravariant deformation functor $F_{\phi, A}'$ (resp.\ $F_{\phi}'$) from $\hat{\C}$ to sets
which associates to $\Omega$ the set of equivalence classes of deformations 
in $F_{\phi, A}(\Omega')$ (resp.\ $F_{\phi}(\Omega')$) where $i: \Omega' \to \Omega$ is any finite purely inseparable morphism.
Likewise, consider the contravariant $A$-Galois functor $F_{A \rtimes_\iota \mu_m, \sigma}'$ from $\hat{\C}$ to sets
which associates to $\Omega$ the set of equivalence classes of covers in
$F_{A \rtimes_\iota \mu_m, \sigma}(\Omega')$ where $i: \Omega' \to \Omega$ is any finite purely inseparable morphism.
The functors $F_\phi'$, $F_{\phi, A}'$ and $F_{A \rtimes_\iota \mu_m, \sigma}'$ are defined over the 
category $\hat{\C}$, but descend to functors on $\C'$ since they associate an equivalence 
class of deformations (or covers) to the equivalence class in $\C'$ of an object $\Omega$ of $\hat{\C}$. 

\begin{corollary} \label{Cisofunc}
The functors $F_{\phi, A}'$ and $F_{A \rtimes_\iota \mu_m, \sigma}'$ are isomorphic. 
\end{corollary}

\begin{proof}
This is automatic since the isomorphism in Theorem \ref{Tfunctor} is compatible with 
finite purely inseparable morphisms.  
\end{proof}

\subsection{A moduli space for $A \rtimes_\iota \mu_m$-Galois covers}  \label{Smoduli}
  
We construct a fine moduli space $\CM_{A \rtimes_\iota \mu_m, \sigma}$ for the functor $F_{A \rtimes_\iota \mu_m, \sigma}'$ 
and give an explicit formula for its Krull dimension.  (The case when $A=\ZZ/p$ appears already in \cite[Thm.\ 2.2.10]{Pr:deg}).
The Krull dimension can be computed using the following formula 
which depends only on numbers which arise from the ramification filtration.

\begin{definition} \label{Ddim}
Given $A \rtimes_\iota \mu_m$ and $\sigma$, let $q=|A|$ and let $s_\iota$ be as in Lemma \ref{Lcong}. 
Let $n(A \rtimes_\iota \mu_m, \sigma)=\#
\{\ell \in \NN^+| \ q \nmid \ell, \ \ell/\gc(\ell,q) \leq m \sigma, \ \ell \equiv s_\iota \bmod m\}$.
\end{definition}

\begin{definition}
Let $n=n(A \rtimes_\iota \mu_m, \sigma)$.
Let $\hat{\GG}_a$ be the formal completion of the group scheme $\GG_a$ at the origin.
Consider the action of $\zeta \in \FF_q^*$ on $(\hat{\GG}_a)^n$ so that $\zeta \circ \vec{r} = \zeta \vec{r}$.  
Let $\CM_{A \rtimes_\iota \mu_m, \sigma}$ be the quotient of $(\hat{\GG}_a)^n$
by this action of $\FF_q^*$.  
\end{definition}

\begin{theorem} \label{Tmoduli}
The functor $F_{A \rtimes_\iota \mu_m, \sigma}'$ is represented by $\CM_{A \rtimes_\iota \mu_m, \sigma}$
which has Krull dimension $n(A \rtimes_\iota \mu_m, \sigma)$.
\end{theorem}

\begin{proof}
To ease notation, let $n=n(A \rtimes_\iota \mu_m, \sigma)$, 
$F'=F_{A \rtimes_\iota \mu_m, \sigma}'$ and $\CM=\CM_{A \rtimes_\iota \mu_m, \sigma}$.
Consider an object $\Omega$ of $\C'$ and the equal characteristic complete local ring $R=\CO_\Omega$.
The following two properties imply that $F'$ is represented by $\CM$.
\begin{description}
\item{i)} There is a natural morphism $T': F'(\circ) \to \Hom_{\C'}(\circ, \CM)$;

\item{ii)} There is a morphism $T: \Hom_{\C'}(\circ, \CM) \to F'(\circ)$ which is an inverse to $T'$.
\end{description}

For $(i)$, suppose we are given $\alpha \in F'(\Omega)$.
Then $\alpha$ corresponds to the isomorphism class of an $A$-Galois cover $\psi_\alpha$ of $X_{\Omega'}$
where $\Omega' \to \Omega$ is a finite purely inseparable morphism.
As in Section \ref{Seq}, $\psi_\alpha$ is given by an equation of the form $v^q-v=r_\alpha$ where $r_\alpha \in R'[[x]][x^{-1}]$
and $R'=\CO_{\Omega'}$.
  
Within the isomorphism class of the cover $\psi_\alpha$, 
there is a choice of equation $v^q-v=r_\alpha'$ with $r_\alpha' \in x^{-1}R'[x^{-1}]$.
This follows from Lemma \ref{Lnoniso} since any element $r \in R'[[x]]$ is of the form $r=d^q-d$ 
where $d \in R'[[x]]$ equals $\sum_{i=0}^{\infty} r^{q^i}$.
   
Within the equivalence class of $\psi_\alpha$, there is a choice of equation
$v^q-v=r_\alpha''$ in standard form.  
Namely, $r_\alpha'' \in x^{-1}R''[x^{-1}]$ and the coefficient of $(x^{-1})^\ell$ in $r''_\alpha$ is zero if $q | i$.
Here $R''$ is a finite purely inseparable extension of $R'$.
The reason is that after a finite radicial extension,  
there is such an $r_\alpha''$ so that $r_\alpha'-r_\alpha''=d_0^{q}-d_0$ for some $d_0 \in R''[x^{-1}]$.  
The pullback of $\psi_\alpha$ to $\Omega''=\Spec(R'')$ is isomorphic to $v^q-v=r_\alpha''$.

Let $s$ be the conductor of $\psi_\alpha$.  By definition, $s \leq m\sigma$.
By Lemma \ref{Lppdegree}, $s$ is the prime-to-$p$ degree of $r_\alpha$, 
One can show this is the same as the prime-to-$p$ degree of $r_\alpha''$.
This implies that the coefficient of $(x^{-1})^\ell$ in $r''_\alpha$ is zero unless $\ell/\gc(\ell,q) \leq m\sigma$.

By Lemma \ref{Lcong}, the coefficient of $(x^{-1})^\ell$ in $r''_\alpha$ is zero unless $\ell \equiv s_\iota$ modulo $m$.
To summarize, the number of exponents $\ell$ so that the coefficient of $(x^{-1})^\ell$ in $r_\alpha''$ can be non-zero 
is the cardinality of
$\{\ell \in \NN^+| \ q \nmid \ell, \ \ell/\gc(\ell,q) \leq m \sigma, \ \ell \equiv s_\iota \bmod m\}$, namely $n$.
The cover $\psi_{\Omega''}$ now yields a morphism $\Omega'' \to (\hat{\GG}_a)^n$ 
using the coefficients of $(x^{-1})^\ell$ from the polynomial $r_\alpha''$.
This yields a morphism $f_\alpha: \Omega'' \to \CM$. 

The morphism $f_\alpha$ is well-defined and uniquely determined by $\alpha$.
First, the image of $f_\alpha$ lies in $\CM$ since the closed fibre of $\psi_\alpha$ is isomorphic to $\psi_{0,\omega}$.
Second, $f_\alpha$ does not depend on the choice of $r_\alpha''$. 
The reason is that if $v_1^q-v_1=r_1(x)$ and $v_2^q-v_2=r_2(x)$ are in standard form,  
then Lemma \ref{Lnoniso} implies that the corresponding $A$-Galois covers $\psi_1$ and $\psi_2$ are isomorphic 
if and only if $r_2(x)=\zeta r_1(x)$ for some $\zeta \in \FF_q^*$.

In the category $\C'$ (but not in the category $\C$), the morphism $f_\alpha$ descends to a morphism 
$\Omega \to \CM$.  Define $T'(\alpha)= f_\alpha \in {\mathop{\rm Hom}}_{\C'}(\Omega, \CM)$.

Conversely, for (ii), suppose
$f \in {\mathop{\rm Hom}}_{\C'}(\Omega, \CM)$.
In other words, suppose $f:\Omega' \to \CM$ where $\Omega'$ is a finite purely inseparable extension of $\Omega$. 
Consider a lifting of $f$ to $(\hat{\GG}_a)^n$. 
Define $r'_f \in x^{-1}R'[x^{-1}]$ using the coordinates of the $\Omega'$-point for the coefficients
of the terms $(x^{-1})^\ell$ for 
$\{\ell \in \NN^+| \ q \nmid \ell, \ \ell/\gc(\ell,q) \leq m \sigma, \ \ell \equiv s_\iota \bmod m\}$.
By analogous arguments as for (i), the polynomial $r'_f$ yields a cover $\psi_{f, \Omega'}$ in standard form.
By Lemma \ref{Lnoniso}, the cover $\psi_{f, \Omega'}$ does not depend on the choice of lifting of $f$.  
Let $T(f)$ to be the equivalence class of $\psi_{f, \Omega'}$ in $F'(\Omega)$.  

The morphisms $T$ and $T'$ are functorial and inverses of each other.
Having verified conditions (i)-(ii), it follows that
$\CM$ represents the functor $F'$ on the category $\C'$.
The moduli space is a scheme of Krull dimension $n$.
\end{proof}

\subsection{Application to \cite{CK}: Inertia $(\ZZ/p)^e \rtimes_\iota \mu_m$} \label{CK}

Suppose $I = (\ZZ/p)^e \rtimes_\iota \mu_m$ and $\phi: Y_k \to U_k$ is an $I$-Galois cover of normal irreducible germs of curves.
This type of inertia group occurs when one studies covers of projective curves which are ordinary.
In Corollary \ref{Creducible}, we use Theorem \ref{Tmoduli} to describe the moduli space parametrizing deformations of $\phi$ 
and give an exact formula for its Krull dimension.
This corollary is a generalization of the dimension count of \cite[Theorem 5.1(a)]{CK}
in which the restriction that the lower jump is $1$ is removed.  

\begin{corollary} \label{Creducible}
Suppose $I =(\ZZ/p)^e \rtimes_\iota \mu_m$.
Let $P_i$ for $1 \leq i \leq r$ be non-trivial elementary abelian $p$-groups 
which are stable and irreducible under the action of $\mu_m$ so that $I=(\times_{i=1}^r P_i) \rtimes_\iota \mu_m$.
Suppose $\phi: Y_k \to U_k$ is an $I$-Galois cover of normal irreducible germs of curves 
with upper jump $\tilde{\sigma}_i$ associated to $g \in P_i$.
In the category $\C'$, the functor $F'_\phi$ is represented by 
$\times_{i=1}^r \CM_{P_i \rtimes_\iota \mu_m, \tilde{\sigma}_i}$
which has Krull dimension $d_\phi = \sum_{i=1}^r n(P_i \rtimes_\iota \mu_m, \tilde{\sigma}_i)$.
\end{corollary}

The formula for $n(P_i \rtimes_\iota \mu_m, \tilde{\sigma}_i)$ is in Definition \ref{Ddim}.

\begin{proof}
Given $I =(\ZZ/p)^e \rtimes_\iota \mu_m$, consider the morphism $\iota: \mu_m \to \Aut((\ZZ/p)^e)$.  
Since $p \nmid m$, Maschke's Theorem implies that this representation is completely reducible.
In other words, there exist non-trivial elementary abelian $p$-groups $P_i$ for $1 \leq i \leq r$ 
which are stable and irreducible under the action of $\mu_m$ so that $I=(\times_{i=1}^r P_i) \rtimes_\iota \mu_m$.

Consider the $P$-Galois subcover $\phi^P:Y \to X$ of $\phi$.
For $1 \leq i \leq r$, let $\phi_i:Y_i \to X$ be the $P_i$-Galois quotient of $\phi^P$.
Then $Y_i$ is invariant under the $\mu_m$-action and 
the cover $\kappa \circ \phi_i$ is $P_i \rtimes_\iota \mu_m$-Galois. 
The ramification filtration of $\kappa \circ \phi_i$ has one jump, which is $\tilde{\sigma}_i$ in the upper numbering.
The subgroup $P_i$ satisfies all the conditions of Lemma \ref{LhypA} for the cover $\kappa \circ \phi_i$.

Consider the functor $F'_{\kappa \circ \phi_i}$ parametrizing equivalence classes of deformations of $\kappa \circ \phi_i$.
The $\mu_m$-Galois quotient of any such deformation must be constant so $F'_{\kappa \circ \phi_i}=F'_{\kappa \circ \phi_i, P_i}$.
The latter functor is isomorphic to $F'_{P_i \rtimes_\iota \mu_m,\tilde{\sigma}_i}$ by Corollary \ref{Cisofunc}. 
By Theorem \ref{Tmoduli}, in the catagory $\C'$, the moduli space representing this last functor 
is $\CM_{P_i \rtimes_\iota \mu_m, \tilde{\sigma}_i}$ which has Krull dimension $n(P_i \rtimes_\iota \mu_m, \tilde{\sigma}_i)$.

The $P$-Galois subcover $\phi^P:Y \to X$ of $\phi$ is the fibre product of the covers $\phi_i:Y_i \to X$.
Giving a deformation of $\phi$ is equivalent to giving a $\mu_m$-invariant deformation of $\phi^P$
which is equivalent to giving deformations of the collection of covers $\kappa \circ \phi_i$ for $1 \leq i \leq r$.    
Thus the moduli space representing the functor $F'_\phi$ is the direct product 
of the moduli spaces representing the functors $F'_{\kappa \circ \phi_i}$. 
\end{proof}

\begin{remark}
In \cite{CK}, Cornelissen and Kato study a functor of deformations of a germ $Y$ with a wildly ramified action 
$\rho:I \to \Aut(Y)$ of the type occuring when $Y$ is the germ of a smooth projective {\it ordinary} curve. 
They compute the formal deformation space of the functor.
In particular, they find its pro-representable hull $H_{\rho}$ and compute its Krull dimension $d_{\rho}$.
They also find the tangent space of the deformation functor and compute its dimension.
The latter dimension may be larger than $d_{\rho}$ 
due to the presence of nilpotent elements coming from obstructed deformations, 
i.e.\ liftings of $(Y, \rho)$ to $k[t]/t^2$ which do not lift to $k[[t]]$.

Apriori, the branch locus might split under deformations of $(Y, \rho)$.
However, it is impossible for the branch locus to split in the ordinary case
since the lower jump of a wildly ramified cover cannot be smaller than $1$.  

The next corollary reproves the dimension count of \cite[Theorem 5.1(a)]{CK}.
A careful analysis shows that the results agree when $d_\rho +1=d_\phi$ (or $d_\rho=d_\phi$ in the case $p =2$ and $e=1$). 
One can check that the results agree except when $p=2$, $m=1$, and $e >1$.  
The calculation of \cite[4.4.3]{CK} does not appear to be accurate when $p=2$, $m=1$, and $e =2$.
\end{remark}

\begin{corollary} \label{Cordinary}
Suppose $\phi$ is the germ of a cover of smooth projective ordinary curves.
Suppose the inertia group $I$ of $\phi$ has order $p^em$ with $p \nmid m$.
Let $c=[\FF_p(\zeta_m), \FF_p]$.  
Then $F'_\phi$ is represented by a scheme of Krull dimension $d_\phi=e/c$.
\end{corollary}

\begin{proof}
The ordinary hypothesis forces the second ramification group in the lower numbering to be trivial, \cite[Thm.\ 2(i)]{Na}.
Thus $I$ is of the form $(\ZZ/p)^e \rtimes_\iota \ZZ/m$ by \cite[IV, Cor.\ 3-4]{Se:lf}. 
The proof now rests only on these two facts and not on the ordinary hypothesis itself.

By Maschke's Theorem, there exist non-trivial elementary abelian $p$-groups $P_i$ for $1 \leq i \leq r$ 
which are stable and irreducible under the action of $\mu_m$ so that $I=(\times_{i=1}^r P_i) \rtimes_\iota \mu_m$.
The upper jump of $\phi$ for $g \in P_i$ is $\tilde{\sigma}_i=1/m$ by Herbrand's formula.
By Corollary \ref{Creducible}, in the category $\C'$, the functor $F'_\phi$ is represented by 
$\times_{i=1}^r \CM_{P_i \rtimes_\iota \mu_m, 1/m}$
which has Krull dimension $d_\phi = \sum_{i=1}^r n(P_i \rtimes_\iota \mu_m, 1/m)$.

Let $\kappa$ be the $\mu_m$-Galois quotient of $\phi$.
The conductor of the $P_i$-subquotient of $\phi$ is $s=1$. 
Thus the integer associated to $\kappa$ and the group $P_i \rtimes_\iota \mu_m$ in Lemma \ref{Lcong} is $s_{\iota,i}=1$.  
Lemma \ref{Lcong}(ii) implies that $|P_i|=p^c$ for $1 \leq i \leq r$ so $e=cr$.

Now $n(P_i \rtimes_\iota \mu_m, 1/m)=
\# \{\ell \in \NN^+| \ p^c \nmid \ell, \ \ell/\gc(\ell,p^c) \leq 1, \ \ell \equiv 1 \bmod m\}$.
The condition $\ell \leq \gc(\ell,p^c)$ implies that $\ell=p^b$ for some $b \in \NN^+$.
Then $p^b \equiv 1 \bmod m$ if and only if $c|b$.
Thus $n(P_i \rtimes_\iota \mu_m, 1/m)=1$ for $1 \leq i \leq r$ and $d_\phi=e/c$.
\end{proof}

\section{Towers of deformations} \label{Stower}

In Section \ref{Sbound}, we give upper and lower bounds for the Krull dimension $d_\phi$ 
of the moduli space parametrizing equiramified deformations of $\phi$ in the category $\C'$.  
These bounds depend only on the ramification filtration of $\phi$.  
We show that the upper bound is realized in the case when $I$ is an abelian $p$-group in Section \ref{Scyclic}.     

The proof involves the study of deformations of towers of covers 
using Corollary \ref{Cisofunc} and Theorem \ref{Tmoduli} along with induction.
The crux issue is whether it is possible to construct a deformation of $\phi$ which dominates a given 
deformation of $\overline{\phi}$ and which is still equiramified.  
We reformulate this issue as an equiramified embedding problem in Section \ref{Sequi}.
The upper bound for $d_\phi$ is realized when $I$ is an abelian $p$-group
since the equiramified embedding problem has a solution when $I$ is a cyclic $p$-group. 
The upper bound for $d_\phi$ is not always realized 
since the equiramified embedding problem does not have a solution in general.

Here is a basic explanation of this issue.  
Every $P$-Galois cover $\phi$ of $X$ is a tower of Artin-Schreier equations $y_i^p-y_i=r_i$.
Every deformation of $\phi$ which maintains the $P$-Galois action 
has the form $y_i^p-y_i=r_i + r_{\alpha,i}$ for some $r_{\alpha,i} \in \CO_\Omega[x^{-1}]$.   
If the deformation is equiramified then the prime-to-$p$ degree of each $r_{\alpha,i}$ in $x^{-1}$ is bounded 
by the corresponding upper jump in the ramification filtration.
The converse is true if $P$ is cyclic, but is false in general. 
Namely, one can not guarantee that the deformation is equiramified
by placing bounds on the degree of each $r_{\alpha,i}$. 
Section \ref{Squaternion} contains an example of this phenomenon for deformations of 
a supersingular elliptic curve in characteristic 2 with an action by the quaternion group.

\subsection{An equiramified embedding problem.} \label{Sequi}

The following problem arises in constructing deformations of towers of covers.

\begin{definition}
An {\it equiramified embedding problem} consists of the following data:\\
$\cdot$ a group $I=P \rtimes_\iota \mu_m$ with $|P|=p^e$ and $p \nmid m$;\\
$\cdot$ a subgroup $A \subset I$ satisfying the conditions of Lemma \ref{LhypA};\\
$\cdot$ an $I$-Galois cover $\phi: Y \to U_k$ of normal irreducible germs of $k$-curves;\\
$\cdot$ a pointed $k$-scheme $(\Omega,\omega)$;\\
$\cdot$ an equiramified deformation $\overline{\phi}_{\Omega}$ of the $I/A$-Galois quotient of $\phi$
over $(\Omega,\omega)$.

A {\it solution} to the equiramified embedding problem consists of  
an equiramified deformation $\phi_{\Omega}$ of $\phi$ over $(\Omega,\omega)$
which dominates $\overline{\phi}_{\Omega}$. 
\end{definition}

\begin{lemma} \label{Lsolution}
Given an equiramified embedding problem over $(\Omega, \omega)$, the subscheme $\Omega'$ of $\Omega$
over which the problem has a solution is non-empty and closed.
\end{lemma}

\begin{proof}
We see that $\Omega'$ is non-empty since it contains $\omega$.
It is closed since the equiramified condition depends only on the generic point of $\Omega'$.
\end{proof}

Every equiramified embedding problem for cyclic $p$-group covers has a solution using class field theory.
We use this in Section \ref{Scyclic} to find the Krull dimension of the moduli space representing equiramified deformations of 
abelian $p$-group covers. 

\begin{lemma} \label{Lcyclicsolution}
Any equiramified embedding problem for which $I$ is a cyclic $p$-group has a solution. 
\end{lemma}

\begin{proof}
By hypothesis, $I \simeq \ZZ/p^e$ for some $e \geq 1$.  
Let $\sigma_1, \ldots, \sigma_e$ be the upper jumps of $\phi$.  
Recall from \cite{Schmid} (see also \cite[Lemma 19]{Pr:genus}) that the upper jumps of a cyclic cover satisfy the 
following condition: either $\sigma_{i+1} =p \sigma_i$ or $\sigma_{i+1} > p \sigma_i$ and $p \nmid \sigma_{i+1}$. 

If $A$ satisfies the conditions of Lemma \ref{LhypA}, 
then $A$ is the last non-trivial higher ramification group $I^{\sigma_e}$ of $\phi$ and $A \simeq \ZZ/p$.
The $I/A$-Galois quotient cover $\overline{\phi}:\overline{Y} \to X_k$ of $\phi$ 
is a $\ZZ/p^{e-1}$-Galois cover of normal irreducible germs of curves. 
The ramification filtration of $\overline{\phi}$ has upper jumps $\sigma_1, \ldots \sigma_{e-1}$.

Consider the given equiramified deformation $\overline{\phi}_{\Omega}$ of $\overline{\phi}$ over $(\Omega, \omega)$.
The fibre of $\overline{\phi}_{\Omega}$ over $\omega$ is isomorphic to $\overline{\phi}$.
By \cite[X, Theorem 5.1]{AGV}, there exists a $\ZZ/p^e$-Galois cover $\phi'_{\Omega}$ dominating $\overline{\phi}_{\Omega}$
whose branch locus is $\xi_\Omega'$.
One can choose $\phi'_{\Omega}$ so that its conductor $s'$ over the generic geometric fibre of $\Omega$
is minimal among all $\ZZ/p^e$-Galois covers dominating $\overline{\phi}_{\Omega}$.
By \cite[Lemma 19]{Pr:genus}, $s' = p \sigma_{e-1}$.  Thus $\sigma_e \geq s'$.

Consider the restriction $\phi'$ of $\phi'_{\Omega}$ over $\omega$.
Then $\phi'$ differs from $\phi$ by an element $\alpha \in \Hom(\pi_1(X'_k), \ZZ/p)$ by \cite[Lemma 3]{Pr:genus}.
Let $s$ be the conductor of $\alpha$.
If $\phi'$ is reducible then $s=\sigma_e$ by \cite[Lemma 9]{Pr:genus}.
If $\phi'$ is irreducible, let $s''$ be the conductor of its normalization.
Then $s'' \leq s'$ since conductors decrease under specialization 
and $s' \leq s''$ by \cite[Lemma 19]{Pr:genus}.  So $s''=s'$.
In this case, $s=\max\{s',\sigma_e\}=\sigma_e$ by \cite[Proposition 7]{Pr:genus}.

Let $\alpha_{0,\Omega}$ be the constant deformation of $\alpha$ over $\Omega$.
Let $\phi_{\Omega}$ be the cover $\phi'_\Omega$ modified by the action of $\alpha_{0,\Omega}$.
The restriction of $\phi_{\Omega}$ over the closed point of $\Omega$ is isomorphic to $\phi$.
By \cite[Proposition 7]{Pr:genus}, $\sigma_e$ is the conductor of $\phi_{\Omega}$ over the generic geometic fibre of $\Omega$.
Thus $\phi_{\Omega}$ is a solution to the equiramified embedding problem.
\end{proof}

Section \ref{Squaternion} contains an example of an equiramified embedding problem 
which does not have a solution for the non-abelian quaternion group of order 8.

\subsection{Equiramified deformations of arbitrary covers} \label{Sbound}

Suppose $I=P \rtimes_\iota \mu_m$ with $|P|=p^e$ and $p \nmid m$.
In this section, we consider equiramified deformations of an $I$-Galois cover $\phi: Y_k \to U_k$ 
of normal irreducible germs of $k$-curves.
We use Corollary \ref{Cisofunc} and Theorem \ref{Tmoduli} to find an upper and lower bound 
for the Krull dimension $d_\phi$ of the moduli space parametrizing deformations of $\phi$ in the category $\C'$. 

The cover $\phi$ can be factored into a tower of covers using its ramification filtration. 
The bottom step of the tower is the $\mu_m$-Galois cover $\kappa:X \to U$.
Suppose $\sigma_1, \ldots, \sigma_e$ are the jumps in the ramification filtration $I^c$
of $\phi$ in the upper numbering. 
The Galois groups of the other steps of the tower are elementary abelian $p$-groups,
namely the quotients $I^{\sigma_{i}}/I^{\sigma_{i+1}}$ of higher ramification groups.
We now modify the indexing of the ramification filtration of $\phi$ to make sure that these quotients
are irreducible under the action of $\mu_m$.
 
Let $1 \subset A_i' \subset A_{i-1}' \ldots \subset A_1'=P$ be the distinct
subgroups of $P$ occuring in the sequence of higher ramification groups of ${\phi}$.
Then $A_i'$ is normal in $I$ so the action of $\mu_m$ stabilizes $A_i'$.
Also $A_i'$ is normal in $A_{i-1}'$ and the quotient is a
non-zero elementary abelian $p$-group.

\begin{definition} \label{Drrf}
The {\it reduced ramification filtration} of ${\phi}$ is a
refinement $1=A_{r+1} \subset A_r \subset A_{r-1} \ldots \subset A_1=P$ of the
sequence of distinct subgroups occuring in the ramification filtration
satisfying the property that $A_{i}/A_{i+1}$ is nontrivial and irreducible under the action of $\mu_m$.
The {\it reduced set of upper jumps} of $\phi$ is the set 
$\{\tilde{\sigma}_i \ | \ 1 \leq i \leq r\}$ 
where $\tilde{\sigma}_i$ is the upper jump associated to $g \in A_{i} - A_{i+1}$.
\end{definition}

The set $\{\tilde{\sigma}_i \ | \ 1 \leq i \leq r\}$ is a subset of 
$\{\sigma_i \ | \ 1 \leq i \leq e\}$.
The multiplicity of the jump $\sigma_i$
in the reduced set of upper jumps is equal to the length of the composition series of
the $\mu_m$-module $I^{\sigma_i}/I^{\sigma_{i+1}}$.

Since $A_{i+1}$ is normal in $I$, the quotient cover $\phi_i$ of $\phi$ by $A_{i+1}$
is $I/A_{i+1}$-Galois for $1 \leq i \leq r$.
Then $A_i/A_{i+1}$ satisfies all the conditions of Lemma \ref{LhypA} with reference to the cover $\phi_i$.  
For example, $A_i/A_{i+1}$ is in the last non-trivial higher ramification group of $\phi_i$.
Let $q_i=|A_i/A_{i+1}|$.

By Lemma \ref{Lcong}, there is a well-defined choice of $s_{\iota,i}$ 
determined from the $\mu_m$-Galois quotient $\kappa$ of $\phi$ and the group $(A_i/A_{i+1}) \rtimes_\iota \mu_m$.
As in Definition \ref{Ddim}, let $n_i=n((A_i/A_{i+1}) \rtimes_\iota \mu_m, \tilde{\sigma}_i)=
\# \{\ell \in \NN^+| \ q_i \nmid \ell, \ \ell/\gc(\ell,q_i) \leq m \tilde{\sigma}_i, \ \ell \equiv s_{\iota,i} \bmod m\}$.

\begin{remark}
By \cite[Lemma 6(i)]{Pr:genus}, the integers $s_{\iota,i}$ 
can also be determined from the ramification filtration of the $I$-Galois cover $\phi$.
For example, if $\phi$ has last lower jump $j_e$, then $s_{\iota,r} \equiv j_e/ |\overline{P}| \bmod m$.   
\end{remark}

\begin{theorem} \label{Tsmoothdef2}
Suppose $\phi: Y_k \to U_k$ is an $I$-Galois cover of normal irreducible germs of curves 
with reduced set of upper jumps $\tilde{\sigma}_1, \ldots \tilde{\sigma}_r$ 
in its reduced ramification filtration $A_r \subset A_{r-1} \ldots \subset A_1$.
In the category $\C'$, the functor $F'_\phi$ is represented by a subscheme $\CM_\phi$ of 
$\times_{i=1}^r \CM_{(A_{i}/A_{i+1}) \rtimes_\iota \mu_m, \tilde{\sigma}_i}$
whose Krull dimension $d_\phi$ satisfies $n_r \leq d_\phi \leq \sum_{i=1}^r n_i$.
\end{theorem}

We expect that an exact formula for $d_\phi$ will depend on the ramification filtration of $\phi$
but not upon $\phi$ itself.
When $m=1$, the upper bound for $d_\phi$ is $\sum_{i=1}^e \lfloor \sigma_i \rfloor - \lfloor \sigma_i/p \rfloor$.

\begin{proof}
The proof is by induction on $r$.  If $r=1$, the result is immediate
from Corollary \ref{Cisofunc} and Theorem \ref{Tmoduli}.

For $r > 1$, let $A=A_r$ be the smallest subgroup in the reduced ramification filtration for $\phi$.
Let $\overline{I}=I/A =P/A \rtimes_\iota \mu_m$.
Consider the $\overline{I}$-Galois quotient $\overline{\phi}=\phi_{r-1}$ of $\phi$.
All information about $\phi$ which is relevant for $\overline{\phi}$ is preserved under this quotient.
Namely, $\overline{\phi}$ has reduced ramification filtration  $A_{r-1}/A_r \subset A_{r-2}/A_r \ldots \subset A_1/A_r$
and reduced set of upper jumps $\tilde{\sigma}_1, \ldots \tilde{\sigma}_{r-1}$.  
By the third isomorphism theorem, $(A_i/A_r)/(A_{i+1}/A_r)=A_i/A_{i+1}$ for $1 \leq i \leq r-1$. 
As a result, the numbers $s_{\iota,i}$ and $n_i$ are the same for $\phi$ and $\overline{\phi}$ for $1 \leq i \leq r-1$.

By the inductive hypothesis, in the category $\C'$, 
the functor $F'_{\overline{\phi}}$ is represented by a subscheme 
$\CM_{\overline{\phi}}$ of $\times_{i=1}^{r-1} \CM_{(A_{i}/A_{i+1})\rtimes_\iota \mu_m, \tilde{\sigma}_i}$
whose Krull dimension $d_{\overline{\phi}}$ satisfies $n_{r-1} \leq d_{\overline{\phi}} \leq \sum_{i=1}^{r-1} n_i$.
Consider the universal $\overline{I}$-Galois deformation 
$\overline{\phi}_{\CM_{\overline{\phi}}}:Y_{\CM_{\overline{\phi}}} \to U_{\CM_{\overline{\phi}}}$
of $\overline{\phi}$ corresponding to this moduli space.

Consider the equiramified embedding problem determined by $\overline{\phi}_{\CM_{\overline{\phi}}}$ and $\phi$.
Let $\overline{\CM}$ be the subscheme of $\CM_{\overline{\phi}}$ over which it has a solution.  
By Lemma \ref{Lsolution}, $\overline{\CM}$ is non-empty and closed and thus 
its irreducible components have Krull dimension between $0$ and $\sum_{i=1}^{r-1} n_i$.
Let $\phi_{\overline{\CM}}$ be a solution to the equiramified embedding problem over $\overline{\CM}$.
The conductor of $\phi_{\overline{\CM}}$ over the generic point of $\overline{\CM}$ is $\tilde{\sigma}_r$ by definition.
The restriction of $\phi_{\overline{\CM}}$ to $\omega$ is isomorphic to $\phi$ by definition.

By Corollary \ref{Cisofunc} and Theorem \ref{Tmoduli}, the functor $F'_{\phi, A}$ is 
represented by $\CM_{A_r \rtimes_\iota \mu_m, \tilde{\sigma}_r}$.
Consider $\CM_\phi=\overline{\CM} \times \CM_{A_r \rtimes_\iota \mu_m, \tilde{\sigma}_r}$.
Then $\CM_\phi$ is a subscheme of $\times_{i=1}^r \CM_{(A_{i}/A_{i+1}) \rtimes_\iota \mu_m, \tilde{\sigma}_i}$
and the Krull dimension of each of its componenets is between $n_r$ and $\sum_{i=1}^r n_i$.
We will show that $\CM_\phi$ represents the functor $F'_\phi$ in the category $\C'$.

First, there is a non-canonical morphism $T': F_\phi'(\circ) \to {\mathop{\rm Hom}}_{\C'}(\circ, \CM_\phi)$.
To see this, let $\Omega$ be an object of $\C'$ and consider $\phi_{\Omega} \in F_\phi'(\Omega)$.
Let $\overline{\phi}_{\Omega}$ be the $\overline{I}$-Galois quotient of $\phi_{\Omega}$.
By the inductive hypothesis, (after inverting by a finite radicial morphism) 
$\overline{\phi}_{\Omega}$ determines a unique morphism $f_1:\Omega \to \CM_{\overline{\phi}}$.
The image of $f_1$ is in $\overline{\CM}$ since $\overline{\phi}_{\Omega}$ can be dominated by the equiramified 
deformation $\phi_{\Omega}$. 

Consider the $I$-Galois cover $f_1^* \phi_{\overline{\CM}}$ of relative $\Omega$-curves.
There exists an element $\alpha \in H_A(\overline{\CM})=\Hom_{\overline{\CM}}(\pi_1(X'_{\overline{\CM}}), A)$    
so that the action of $\alpha$ takes $f_1^* \phi_{\overline{\CM}}$ to $\phi_{\Omega}$ by \cite[Lemma 3]{Pr:genus}.
Then $\alpha \in H_A^\iota(\overline{\CM})$ by \cite[Lemma 5]{Pr:genus}.
Then $\alpha \in H_A^{\iota, \tilde{\sigma}_r}$ since the conductor of $\alpha$ is at most $m\tilde{\sigma}_r$
by \cite[Proposition 7]{Pr:genus}.
Finally, $\alpha \in F_{\phi,A}'(\overline{\CM})$ because $f_1^* \phi_{\overline{\CM}}$ and $\phi_{\Omega}$ are both 
isomorphic to $\phi$ over $\omega$.  

By Theorem \ref{Tmoduli}, 
$\alpha$ determines a unique morphism $f_2:\Omega \to \CM_{A_r \rtimes_\iota \mu_m, \tilde{\sigma}_r}$
in the category $\C'$.
Let $f=(f_1, f_2)$.  
The conclusion is that, in the category $\C'$, there is a morphism $f:\Omega \to \CM_\phi$.  Let $T'(\phi_\Omega)=f$.
The morphism $T'$ is non-canonical since it depends on the choice of $\phi_{\overline{\CM}}$.

Secondly, there is a morphism $T:  {\mathop{\rm Hom}}_{\C'}(\circ, \CM_\phi) \to F_\phi'(\circ)$ which is an inverse to $T'$.
To see this, let $f: \Omega \to \CM_\phi$ be an element of ${\mathop{\rm Hom}}_{\C'}(\circ, \CM_\phi)$.  
This yields morphisms $f_1: \Omega \to \overline{\CM}$ and 
$f_2: \Omega \to \CM_{A_r \rtimes_\iota \mu_m, \tilde{\sigma}_r}$.
Then $f_1^* \phi_{\overline{\CM}} \in F_\phi'(\Omega)$ is an equiramified deformation of $\phi$ over $\Omega$. 
In particular, $f_1^* \phi_{\overline{\CM}}$ is an $I$-Galois cover of $U_\Omega$ whose fibre over 
$\omega$ is isomorphic to $\phi$.  
By Corollary \ref{Cisofunc} and Theorem \ref{Tmoduli}, $f_2$ determines an element $\alpha_\Omega$ of  
$F_{A \rtimes_\iota \mu_m, \sigma}'(\Omega)=F_{\phi,A}'(\Omega)$.  
Define $T(f)$ to be the element of $F_\phi'(\Omega)$ 
obtained by the action of $\alpha_\Omega$ on $f_1^* \phi_{\overline{\CM}}$.  

The morphisms $T$ and $T'$ are functorial.
It follows from the inductive hypothesis and Theorem \ref{Tmoduli} that $T$ and $T'$ are inverses.
So $\CM_\phi$ represents the functor $F_\phi'$ of 
equivalence classes of equiramified deformations of $\phi$ on the category $\C'$.  
\end{proof}

\begin{remark} \label{Rothercompare}
In \cite{Ha:mod}, Harbater constructs a moduli space $M_P$ for $P$-Galois covers of germs of curves.
The moduli space in this paper is different in several respects.
First, the inertia group $I$ in this paper can be a cyclic-by-$p$ group (not just a $p$-group). 
Second, in \cite{Ha:mod} the author uses an equivalence only under \'etale pullbacks so that 
$M_P$ can be defined over the category of schemes. 
Third, in \cite{Ha:mod} there is no equivalence between covers under inseparable pullbacks, 
so $M_P$ is highly singular.    
Finally, in \cite{Ha:mod}, the author considers deformations of \'etale covers of punctured curves.
With these deformations, 
there is no bound for the jumps in the ramification filtration on the generic fibre.   
As a result, $M_P$ is infinite-dimensional.
In contrast, the moduli space in this paper is finite-dimensional.

A versal deformation space for {\it non-Galois} wildly ramified covers of germs of curves is constructed by Fried and Mezard 
in \cite{FM}.  

In \cite{Kont}, Kontogeorgis compares deformations of wildly ramified covers of germs of curves 
to deformations of Galois representations.
\end{remark}

\subsection{Equiramified deformations of abelian $p$-group covers} \label{Scyclic}

Suppose $I$ is an abelian $p$-group.
Suppose $\phi: Y_k \to U_k$ is an $I$-Galois cover of normal irreducible germs of $k$-curves  
with upper jumps $\sigma_1, \ldots \sigma_e$ in its ramification filtration.
We describe the moduli space parametrizing deformations of $\phi$ and give an exact formula for its Krull dimension
in terms of the upper jumps.
The case when $I$ is an abelian $p$-group is simpler than the general case due to class field theory.

\begin{corollary} \label{Tcyclic}
Suppose $I$ is an abelian $p$-group.
Suppose $\phi: Y \to U_k$ is an $I$-Galois cover with upper jumps $\sigma_1, \ldots, \sigma_e$.
In the category $\C'$, the functor $F'_\phi$ is represented by $\times_{i=1}^e \CM_{\ZZ/p, \sigma_i}$
which has Krull dimension $d_\phi=\sum_{i=1}^e (\sigma_i - \lfloor \sigma_i/p \rfloor)$. 
\end{corollary}

\begin{proof}
Since an abelian $p$-group is the direct product of cyclic $p$-groups,  
the cover $\phi$ is a fibre product of a collection of covers whose inertia groups are cyclic $p$-groups. 
The functor of (equiramified) deformations of $\phi$ is equivalent to the product of the functors of 
(equiramified) deformations of this collection of cyclic quotients. 
The ramification filtration and its jumps in the upper numbering are preserved under quotients.  
The proof thus reduces to the case that $I \simeq \ZZ/p^e$.

The cyclic case uses induction on $e$; the case $e=1$ is covered by Theorem \ref{Tmoduli}.

Let $A$ be the last non-trivial higher ramification group $I^{\sigma_e}$ of $\phi$.
Since $\sigma_e \geq p \sigma_{e-1}$, 
it follows that $A \simeq \ZZ/p$ and that $A$ satisfies the conditions of Lemma \ref{LhypA}. 
The $I/A$-Galois cover $\overline{\phi}:\overline{Y} \to U_k$  
is a $\ZZ/p^{e-1}$-Galois cover of normal irreducible germs of curves. 
Then $\overline{\phi}$ has upper jumps $\sigma_1, \ldots \sigma_{e-1}$ in its ramification filtration.

Applying the inductive hypothesis to $\overline{\phi}$, in the category $\C'$, 
the functor $F'_{\overline{\phi}}$ is represented by $\times_{i=1}^{e-1} \CM_{\ZZ/p, \sigma_i}$
which has Krull dimension $d_{\overline{\phi}}=\sum_{i=1}^{e-1} n(\ZZ/p, \sigma_i)$. 
The cover $\phi$ and the universal $\overline{I}$-Galois deformation of $\overline{\phi}$ corresponding to this moduli space
determine an equiramified embedding problem.
By Lemma \ref{Lcyclicsolution}, this equiramified embedding problem has a solution 
over $\times_{i=1}^{e-1} \CM_{\ZZ/p, \sigma_i}$.
By Theorem \ref{Tmoduli}, the functor of deformations of $\phi$ 
with constant $I/A$-quotient deformation is represented by $\CM_{\ZZ/p, \sigma_e}$ in the category $\C'$.
It follows from the proof of Theorem \ref{Tsmoothdef2} 
that $F'_\phi$ is represented by $\times_{i=1}^e \CM_{\ZZ/p, \sigma_i}$ in the category $\C'$.

The Krull dimension $d_\phi$ of $\times_{i=1}^e \CM_{\ZZ/p, \sigma_i}$ is $\sum_{i=1}^e n(\ZZ/p, \sigma_i)$.
By Definition \ref{Ddim}, $n(\ZZ/p, \sigma_i)=\# \{\ell \in \NN^+|\  \ell \leq \sigma_i, \ p \nmid \ell\}$.
The Hasse-Arf Theorem implies that $\sigma_i \in \NN$.
Thus $d_\phi=\sum_{i=1}^e (\sigma_i - \lfloor \sigma_i/p \rfloor)$.
\end{proof}

\begin{example} {\bf A minimal $(\ZZ/p^e)$-Galois action.}

There exists a $\ZZ/p^e$-Galois cover $\phi:Y_k \to U_k$ of normal irreducible germs of $k$-curves 
whose ramification filtration has upper jumps $\{1,p,p^2, \ldots, p^{e-1}\}$.
This is the smallest possible sequence of upper jumps for a $\ZZ/p^e$-Galois cover.
By Theorem \ref{Tcyclic}, the Krull dimension
of the moduli space parametrizing equiramified deformations of $\phi$ is 
$d_\phi = p^{e-1}$.
This is because $d_\phi=\sum_{i=1}^e n(\ZZ/p, p^{i-1}) =\sum_{i=1}^e (p^{i-1}-p^{i-2})$.
\end{example}

\begin{remark}
In \cite[Proposition 4.1.1]{BM}, the authors consider deformations of a $\ZZ/p^e$-Galois action on a germ of a curve
over a local Artinian $k$-algebra of mixed characteristic.  
They show that the dimension of the tangent space is $\lfloor 2\beta/p^e \rfloor - \lceil \beta/p^e \rceil$
where $\beta$ is the degree of the discrimininant.
Even for $\ZZ/p$-covers, this shows that not every deformation in equal characteristic lifts to 
a deformation in mixed characteristic.
\end{remark}

\subsection{Future application: the dimension of $\CM_g[G] \otimes \FF_p$} \label{BM}

Suppose $\varphi:\Y \to \X$ is a (wildly ramified) $G$-Galois cover of smooth projective curves branched at $\B$. 
If $b \in \B$, let $d_b$ be the Krull dimension of the moduli space
parametrizing deformations of the germ $\phi_b$ of $\phi$ above $b$.
Theorem \ref{Tsmoothdef2} gives an upper and lower bound for $d_b$ in terms of 
the reduced set of upper jumps of the reduced ramification filtration for $\phi_b$.
Using formal patching, one can show that $\sum_{b\in \B} d_b$ is the Krull dimension of the component 
of the Hurwitz space for wildly ramified covers of $\X$ branched at $\B$ which contains $\varphi$. 

Let $\rho:G \to \Aut(\Y)$ be the action associated to $\phi$. 
Consider the component of the moduli space $\CM_{g}[G] \otimes \FF_p$
of curves of genus $g$ with an action by $G$ which contains $(\Y, \rho)$.
The dimension $\delta_{\Y, \rho}$ of this component is larger than $\sum_{b\in \B} d_b$ for several reasons.  
When studying $\CM_{g}[G] \otimes \FF_p$, one must consider deformations of $(\Y, \rho)$ 
for which the base curve $\X$, the branch locus $\B$, and even the size of $\B$ do not remain constant.  
In particular, the branch locus will be a relative Cartier divisor of fixed degree 
on the relative curve $\X_\Omega$.  

At this time, it appears that the dimension of $\CM_{g}[G] \otimes \FF_p$ is known only 
in the case that $G=\ZZ/p$.
The computation of $\delta_{\Y, \rho}$ in this case relies on the corresponding computation of the dimension $d_b$.
For example, in \cite[Section 5.2]{BM}, the authors consider $(\Y, \rho)$ so that the associated cover 
$\varphi$ is a $\ZZ/p$-Galois cover of the projective line branched at only one point $b$ 
with inertia $\ZZ/p$ and lower jump $j$.  Here $\Y$ has genus $g=(p-1)(j-1)/2$.    
They investigate deformations of $(\Y, \rho)$ for which the one branch point $b$
splits into $h$ branch points $b_1, \ldots, b_h$ with lower jumps $j_{b_1}, \ldots, j_{b_h}$.
They use the dimension $d_{b_i}=j_{b_i}-\lfloor j_{b_i} \rfloor$ to show that the component of 
$\CM_{g}[\ZZ/p] \otimes \FF_p$ containing $(\Y, \rho)$ has dimension $j-2$.   

For this reason, we expect that the formula for $d_b$ from Theorem \ref{Tcyclic}
can be used to determine the dimension of the component of $\CM_g[G] \otimes \FF_p$ containing $(\Y, \rho)$ 
when $G=\ZZ/p^e$ and that an analogous formula will be necessary to compute this dimension in the general case.

\subsection{Example: A minimal quaternion-Galois action.} \label{Squaternion}

We consider deformations of a supersingular elliptic curve with quaternion group action.
Let $k$ be an algebraically closed field of characteristic $p=2$.
Let $\zeta_3$ be a root of $x^2+x+1$ in $k$.
Let $D_8=\{1, \mu, \tau, \mu \tau, [-1], [-1]\mu, [-1] \tau, [-1] \mu \tau\}$ 
be the quaternion group of order 8.
We consider deformations of a $D_8$-Galois cover of the 
projective line branched at only one point with minimal ramification filtration.  

\begin{lemma} \label{Lquat}
The cover $\Phi:Y \to {\mathbb P}^1_k$
given by the following equations has Galois group $D_8$: 
$v^2-v=u$; $w^2+w=v$; $y^2-y=w^3$.
\end{lemma}

\begin{proof}
Define the action of $\mu$ by 
$\mu(v)=v$, $\mu(w)=w+1$, $\mu(y)=y+w+\zeta_3$.
Define the action of $\tau$ by 
$\tau(v)=v+1$, $\tau(w)=w+\zeta_3$, $\tau(y)=y+w(\zeta_3+1)+\zeta_3$.
One can check that $[-1]=\mu^2=\tau^2$ fixes $v$ and $w$, and $[-1](y)=y+1$, and $\mu \tau = [-1] \tau \mu$. 
\end{proof}

\begin{lemma}
The cover $\Phi$ is branched at only one point, where it is totally
ramified.  The jumps in the ramification filtration in the lower
numbering are $1,1,3$ and in the upper numbering are $1,1,3/2$.  The
curve $Y$ is a supersingular elliptic curve.  
\end{lemma}

\begin{proof}
By the Jacobian criterion, the cover $\Phi$ is \'etale over $k[t]$.
Since $D_8$ is a $2$-group, $\Phi$ must be totally ramified over $\infty$.
The lower jumps are the absolute values of the valuations of $u$, $v$, and $w^3$
in respectively $k[u^{-1}]$, $k[v^{-1}]$, and $k[w^{-1}]$.  The calculation of the
upper jumps follows from Herbrand's formula.  The curve $Y$ has genus
$1$ by the Riemann-Hurwitz formula, \cite{Har} \cite[IV, Prop.\ 4]{Se:lf}, 
and is supersingular by the During-Shafarevic formula, \cite{Cr}.     
\end{proof}

Recall that there is a unique elliptic curve which is a $D_8$-Galois cover of
the projective line and that this cover can be modified by a two-dimensional family of 
affine linear transformations.  The resulting two-dimensional family of covers is isotrivial.
We now describe these deformations explicitly. 

\begin{proposition} \label{Pquat}
The equations $v^2-v=u +a_1u$; $w^2+w=v+a_2u$; $y^2-y=w^3+a_3u$
yield a family of $D_8$-Galois covers of the projective line
specializing to $\Phi$ at $(a_1,a_2,a_3)=(0,0,0)$.
The normalization of a fibre of this family yields a smooth connected curve  
as long as $a_1 \not = -1$ and $a_2/(a_1+1)$ is not a root of $x^2+x+1$.
This curve has genus 1 if $a_2=0$ or $a_2=a_1+1$ and genus 2 otherwise.
\end{proposition}

\begin{proof}
One can define the action of $\mu$ and $\tau$ on the variables
$v$,$w$, and $y$ exactly as in the proof of Lemma \ref{Lquat}.  
It follows that this is a $D_8$-Galois cover of relative curves over $k[a_1,a_2,a_3]$
which specializes to $\Phi$ when $(a_1,a_2,a_3)=(0,0,0)$. 

Consider the normalization of a fibre of this family over some
geometric point $(a_1,a_2,a_3)$.  
We denote this cover $Y \to W \to V \to {\mathbb P}^1$ where each of the three
steps in this tower is a degree two Artin-Schreier extension branched at exactly one point.  
We will investigate each of the curves $V$, $W$ and $Y$
to see whether $Y$ is connected.
 
The cover $V \to {\mathbb P}^1$ is given by the equation
$v^2-v=(a_1+1)u$.  This equation is in standard form and one sees
that $V$ is disconnected if and only if $a_1=-1$.  
Let us suppose that $V$ is connected, in which case it is a smooth curve of genus $0$.
Also, $v^{-1}$ is a uniformizer of $V$ at its ramification point.

The cover $W \to V$ is given by 
$w^2+w=v+a_2u=v+(v^2-v)a_2/(a_1+1)$, which is not in standard form
since the exponents of $v$ are not all odd.
To determine whether $W$ is connected
we will change this equation into standard form.  
(This destroys the $D_8$-Galois structure of the tower.)
After a purely inseparable extension, there exists an element $c_1$
so that $c_1^2=a_2/(a_1+1)$.  Let $c_2=1+c_1+c_1^2$.
The ${\mathbb Z}/2$-cover $W \to V$ is isomorphic to the ${\mathbb Z}/2$-cover
$W_1 \to V$ with equation $w_1^2+w_1=vc_2$.
The isomorphism identifies $w$ with $w_1+c_1v$.
It follows that $W$ is disconnected if and only if $c_2=0$ 
which is equivalent to the condition that $a_2/(a_1+1)$ is a root of $x^2+x+1$.
Let us suppose that $W$ is connected, in which case it is a smooth curve of genus $0$.

We now consider the cover $Y \to W \simeq W_1$.
The element $w_1^{-1}$ is a uniformizer of $W_1$ at its ramification point.
One can check that $v=(w_1^2-w_1)/c_2$.  
This implies that $w=c_3w_1^2+c_4w_1$ where $c_3=c_1/c_2$ and $c_4=1-c_1/c_2$.  
We note that $c_3$ and $c_4$ are defined under the restrictions on $a_1$
and $a_2$ and are not simultaneously zero. 
One can also check that $a_3u$ has degree 4 in the variable $w_1$. 

The equation $y^2-y=w^3+a_3u$ for the cover $Y \to W$
yields an equation $y^2-y=(c_3w_1^2+c_4w_1)^3 +a_3u$ for the cover $Y \to W_1$.
After putting this equation in standard form, its leading terms are 
$(c_3^2c_4)w_1^5+(c_4^3+c_3^{3/2})w_1^3$.  The fact that this equation is non-constant 
implies that $Y$ is connected whenever $W$ is.  

The lower jump of $Y \to U$ equals the degree of this equation; in other words, 
it is 3 if $c_3c_4=0$ and is 5 if $c_3c_4 \not =0$.
We see that $c_3=0$ if and only if $a_2=0$ and $c_4=0$ if and only if $a_2=a_1+1$.
It follows that $Y$ is a smooth curve of genus 1 if $a_2=0$ or $a_2=a_1+1$ and of genus 2 otherwise.
\end{proof}

\begin{corollary} \label{Cquat}
Let $R$ be an equal characteristic complete local ring.
Every equiramified deformation of $\Phi$ over $R$ is isomorphic to a fibre of the family 
$v^2-v=u+a_1u$; $w^2+w=v$; $y^2-y=w^3+a_3u$ where $a_1$ and $a_3$ are in the maximal ideal of $R$. 
\end{corollary}

\begin{proof}
By Theorem \ref{Tsmoothdef2}, an equiramified deformation corresponds to a fibre of the family in Proposition \ref{Pquat}
for which $a_1, a_2, a_3$ are in the maximal ideal of $R$.
Then $a_2 \not = a_1+1$ so $a_2=0$.
\end{proof}

This corollary is relevant to the study of deformations of germs of curves (rather than projective curves)
by the theorem of Katz-Gabber \cite{KG}.  
Let $\hat{\Phi}$ be the germ of the cover $\Phi$ above $u=\infty$.
By Theorem \ref{Tsmoothdef2}, $n_3 \leq d_{\hat{\Phi}} \leq n_1+n_2+n_3$.  
By Definition \ref{Ddim}, $n_1=n_2=\#\{\ell \in {\mathbb N}^+ | \ \ell \leq 1, \ 2 \nmid \ell\}=1$ and 
$n_3=\#\{\ell \in {\mathbb N}^+| \  \ell \leq 3/2, \ 2 \nmid \ell\}=1$.
So $1 \leq d_{\hat{\Phi}} \leq 3$.  
In fact, one sees that $d_{\hat{\Phi}}=2$ using the two-dimensional family of deformations 
of $\Phi$ from Corollary \ref{Cquat}.

\bibliographystyle{abbrv} 
\bibliography{paper}

\noindent
Rachel J. Pries\\
Department of Mathematics\\
Colorado State University\\
Fort Collins, CO 80524\\
pries@math.colostate.edu

\end{document}